\documentstyle[12pt,amssymb,amscd]{article}

\newenvironment{spequation}[1]%
{\renewcommand{\theequation}{#1}  
\begin{equation}}  
{\end{equation}%
\renewcommand{\theequation}{(\arabic{equation})}%
\addtocounter{equation}{-1}}     
  
\begin{document}

\baselineskip 22pt     
     
\newtheorem{definition}{\normalsize\sc Definition}[section]     
\newtheorem{prop}[definition]{\normalsize\sc Proposition}     
\newtheorem{lem}[definition]{\normalsize\sc Lemma}     
\newtheorem{corollary}[definition]{\normalsize\sc Corollary}     
\newtheorem{theorem}[definition]{\normalsize\sc Theorem}     
\newtheorem{example}[definition]{\normalsize\sc  Example}     
\newtheorem{remark}[definition]{\normalsize\sc Remark}     
     
\newcommand{\nc}[2]{\newcommand{#1}{#2}}     
\newcommand{\rnc}[2]{\renewcommand{#1}{#2}}     
     
\nc{\Section}{\setcounter{definition}{0}\section}     
\renewcommand{\theequation}{\thesection.\arabic{equation}}     
\newcounter{c}     
\renewcommand{\[}{\setcounter{c}{1}$$}     
\newcommand{\etyk}[1]{\vspace{-7.4mm}$$\begin{equation}\Label{#1}     
\addtocounter{c}{1}}     
\renewcommand{\]}{\ifnum \value{c}=1 $$\else \end{equation}\fi}     
\renewcommand\refname{\normalsize\sc\centering References}

   
\nc{\bpr}{\begin{prop}}     
\nc{\bth}{\begin{theorem}}     
\nc{\ble}{\begin{lem}}     
\nc{\bco}{\begin{corollary}}     
\nc{\bre}{\begin{remark}}     
\nc{\bex}{\begin{example}}     
\nc{\bde}{\begin{definition}}     
\nc{\ede}{\end{definition}}     
\nc{\epr}{\end{prop}}     
\nc{\ethe}{\end{theorem}}     
\nc{\ele}{\end{lem}}     
\nc{\eco}{\end{corollary}}     
\nc{\ere}{\hfill\mbox{$\Diamond$}\end{remark}}     
\nc{\eex}{\end{example}}     
\nc{\epf}{\hfill\mbox{$\Box$}}     
\nc{\ot}{\otimes}     
\nc{\bsb}{\begin{Sb}}     
\nc{\esb}{\end{Sb}}     
\nc{\ct}{\mbox{${\cal T}$}}     
\nc{\ctb}{\mbox{${\cal T}\sb B$}}     
\nc{\bcd}{\[\begin{CD}}     
\nc{\ecd}{\end{CD}\]}     
\nc{\ba}{\begin{array}}     
\nc{\ea}{\end{array}}     
\nc{\bea}{\begin{eqnarray}}     
\nc{\eea}{\end{eqnarray}}     
\nc{\be}{\begin{enumerate}}     
\nc{\ee}{\end{enumerate}}     
\nc{\beq}{\begin{equation}}     
\nc{\eeq}{\end{equation}}     
\nc{\bi}{\begin{itemize}}     
\nc{\ei}{\end{itemize}}     
\nc{\kr}{\mbox{Ker}}     
\nc{\te}{\!\ot\!}     
\nc{\pf}{\mbox{$P\!\sb F$}}     
\nc{\pn}{\mbox{$P\!\sb\nu$}}     
\nc{\bmlp}{\mbox{\boldmath$\left(\right.$}}     
\nc{\bmrp}{\mbox{\boldmath$\left.\right)$}}     
\rnc{\phi}{\mbox{$\varphi$}}     
\nc{\LAblp}{\mbox{\LARGE\boldmath$($}}     
\nc{\LAbrp}{\mbox{\LARGE\boldmath$)$}}     
\nc{\Lblp}{\mbox{\Large\boldmath$($}}     
\nc{\Lbrp}{\mbox{\Large\boldmath$)$}}     
\nc{\lblp}{\mbox{\large\boldmath$($}}     
\nc{\lbrp}{\mbox{\large\boldmath$)$}}     
\nc{\blp}{\mbox{\boldmath$($}}     
\nc{\brp}{\mbox{\boldmath$)$}}     
\nc{\LAlp}{\mbox{\LARGE $($}}     
\nc{\LArp}{\mbox{\LARGE $)$}}     
\nc{\Llp}{\mbox{\Large $($}}     
\nc{\Lrp}{\mbox{\Large $)$}}     
\nc{\llp}{\mbox{\large $($}}     
\nc{\lrp}{\mbox{\large $)$}}     
\nc{\lbc}{\mbox{\Large\boldmath$,$}}     
\nc{\lc}{\mbox{\Large$,$}}     
\nc{\Lall}{\mbox{\Large$\forall$}}     
\nc{\bc}{\mbox{\boldmath$,$}}     
\rnc{\epsilon}{\varepsilon}     
\nc{\ra}{\rightarrow}     
\nc{\ci}{\circ}     
\nc{\cc}{\!\ci\!}     
\nc{\T}{\mbox{\sf T}}     
\nc{\can}{\mbox{\em\sf T}\!\sb R}     
\nc{\cnl}{$\mbox{\sf T}\!\sb R$}     
\nc{\lra}{\longrightarrow}     
\nc{\M}{\mbox{Map}}     
\nc{\imp}{\Rightarrow}     
\rnc{\iff}{\Leftrightarrow}     
\nc{\bmq}{\cite{bmq}}     
\nc{\ob}{\mbox{$\Omega\sp{1}\! (\! B)$}}     
\nc{\op}{\mbox{$\Omega\sp{1}\! (\! P)$}}     
\nc{\oa}{\mbox{$\Omega\sp{1}\! (\! A)$}}     
\nc{\inc}{\mbox{$\,\subseteq\;$}}     
\nc{\de}{\mbox{$\Delta$}}     
\nc{\spp}{\mbox{${\cal S}{\cal P}(P)$}}     
\nc{\dr}{\mbox{$\Delta_{R}$}}     
\nc{\dsr}{\mbox{$\Delta_{\cal R}$}}     
\nc{\m}{\mbox{m}}     
\nc{\0}{\sb{(0)}}     
\nc{\1}{\sb{(1)}}     
\nc{\2}{\sb{(2)}}     
\nc{\3}{\sb{(3)}}     
\nc{\4}{\sb{(4)}}     
\nc{\5}{\sb{(5)}}     
\nc{\6}{\sb{(6)}}     
\nc{\7}{\sb{(7)}}     
\nc{\hsp}{\hspace*}     
\nc{\nin}{\mbox{$n\in\{ 0\}\!\cup\!{\Bbb N}$}}     
\nc{\al}{\mbox{$\alpha$}}     
\nc{\bet}{\mbox{$\beta$}}     
\nc{\ha}{\mbox{$\alpha$}}     
\nc{\hb}{\mbox{$\beta$}}     
\nc{\hg}{\mbox{$\gamma$}}     
\nc{\hd}{\mbox{$\delta$}}     
\nc{\he}{\mbox{$\varepsilon$}}     
\nc{\hz}{\mbox{$\zeta$}}     
\nc{\hs}{\mbox{$\sigma$}}     
\nc{\hk}{\mbox{$\kappa$}}     
\nc{\hm}{\mbox{$\mu$}}     
\nc{\hn}{\mbox{$\nu$}}     
\nc{\la}{\mbox{$\lambda$}}     
\nc{\hl}{\mbox{$\lambda$}}     
\nc{\hG}{\mbox{$\Gamma$}}     
\nc{\hD}{\mbox{$\Delta$}}     
\nc{\ho}{\mbox{$\omega$}}     
\nc{\hO}{\mbox{$\Omega$}}     
\nc{\hp}{\mbox{$\pi$}}     
\nc{\hP}{\mbox{$\Pi$}}     
\nc{\bpf}{{\it Proof.~~}}     
\nc{\slq}{\mbox{$A(SL\sb q(2))$}}     
\nc{\fr}{\mbox{$Fr\llp A(SL(2,\IC))\lrp$}}     
\nc{\slc}{\mbox{$A(SL(2,\IC))$}}     
\nc{\af}{\mbox{$A(F)$}}     
\rnc{\widetilde}{\tilde}     
\nc{\qdt}{quantum double torus}     
\nc{\aqdt}{\mbox{$A(DT^2_q)$}}     
\nc{\dtq}{\mbox{$DT^2_q$}}     
\nc{\uc}{\mbox{$U(2)$}}     
\nc{\uq}{\mbox{$U_{q^{-1},q}(2)$}}     
\rnc{\subset}{\inc}

\def\esl{{\mbox{$E\sb{\frak s\frak l (2,{\Bbb C})}$}}}     
\def\esu{{\mbox{$E\sb{\frak s\frak u(2)}$}}}     
\def\zf{{\mbox{${\Bbb Z}\sb 4$}}}     
\def\zt{{\mbox{$2{\Bbb Z}\sb 2$}}}     
\def\ox{{\mbox{$\Omega\sp 1\sb{\frak M}X$}}}     
\def\oxh{{\mbox{$\Omega\sp 1\sb{\frak M-hor}X$}}}     
\def\oxs{{\mbox{$\Omega\sp 1\sb{\frak M-shor}X$}}}     
\def\Fr{\mbox{Fr}}     
\def\gal{-Galois extension}     
\def\hge{Hopf-Galois extension}     
\def\cge{coalgebra-Galois extension}     
\def\pge{$\psi$-Galois extension}     
\def\ta{\tilde a}     
\def\tb{\tilde b}     
\def\tc{\tilde c}     
\def\td{\tilde d}     
\def\st{\stackrel}     
     
   
\newcommand{\Sp}{{\rm Sp}\,}     
\newcommand{\Mor}{\mbox{$\rm Mor$}}     
\newcommand{\skrA}{{\cal A}}     
\newcommand{\Phase}{\mbox{$\rm Phase\,$}}     
\newcommand{\id}{{\rm id}}     
\newcommand{\diag}{{\rm diag}}     
\newcommand{\inv}{{\rm inv}}     
\newcommand{\ad}{{\rm ad}}     
\newcommand{\poi}{{\rm pt}}     
\newcommand{\Dim}{{\rm dim}\,}     
\newcommand{\Ker}{{\rm ker}\,}     
\newcommand{\Mat}{{\rm Mat}\,}     
\newcommand{\Rep}{{\rm Rep}\,}     
\newcommand{\Fun}{{\rm Fun}\,}     
\newcommand{\Tr}{{\rm Tr}\,}     
\newcommand{\supp}{\mbox{$\rm supp$}}     
\newcommand{\half}{\frac{1}{2}}     
\newcommand{\skrF}{{A}}     
\newcommand{\skrD}{{\cal D}}     
\newcommand{\skrC}{{\cal C}}     
\newcommand{\ttimes}{\mbox{$\hspace{.5mm}\bigcirc\hspace{-4.9mm}     
\perp\hspace{1mm}$}}     
\newcommand{\Ttimes}{\mbox{$\hspace{.5mm}\bigcirc\hspace{-3.7mm}     
\raisebox{-.7mm}{$\top$}\hspace{1mm}$}}     
\newcommand{\bbr}{{\bf R}}     
\newcommand{\bbz}{{\bf Z}}     
\newcommand{\Ci}{C_{\infty}}     
\newcommand{\Cb}{C_{b}}     
\newcommand{\fa}{\forall}     
\newcommand{\rrr}{right regular representation}     
\newcommand{\wrt}{with respect to}     
\newcommand{\qg}{quantum group}     
\newcommand{\qgs}{quantum groups}     
\newcommand{\cs}{classical space}     
\newcommand{\qs}{quantum space}     
\newcommand{\po}{Pontryagin}     
\newcommand{\ch}{character}     
\newcommand{\chs}{characters}     
     
\def\inbar{\,\vrule height1.5ex width.4pt depth0pt}     
\def\IC{{\Bbb C}}     
\def\IZ{{\Bbb Z}}     
\def\IN{{\Bbb N}}     
\def\otc{\otimes_{\IC}}     
\def\ra{\rightarrow}     
\def\ota{\otimes_ A}     
\def\otza{\otimes_{ Z(A)}}     
\def\otc{\otimes_{\IC}}     
\def\h{\rho}     
\def\x{\zeta}     
\def\th{\theta}     
\def\s{\sigma}     
\def\t{\tau}     
     
\def\sw#1{{\sb{(#1)}}}     
\def\sua#1{{\sp{[#1]}}}
\def\sco#1{{\sp{(\bar #1)}}} 
\def\su#1{{\sp{(#1)}}} 
\def\extd{{\rm d}}     
\def\covD{{\rm D}}  
\def\proof{\noindent{\sl Proof.~~}}     
\def\endproof{\hbox{$\sqcup$}\llap{\hbox{$\sqcap$}}\medskip}     
\def\tens{\mathop{\otimes}}     
\def\CC{{\frak C}}     
\def\CL{{\Lambda}}     
\def\M{{\bf M}}     
\def\CM{{\bf  Mod}}  
\def\o{{}_{(1)}}     
\def\t{{}_{(2)}}     
\def\Bo{{}_{\und{(1)}}}     
\def\Bt{{}_{\und{(2)}}}     
\def\th{{}_{(3)}}     
\def\Bth{{}_{\und{(3)}}}     
\def\<{{\langle}}     
\def\>{{\rangle}}     
\def\und#1{{\underline{#1}}}     
\def\ra{{\triangleleft}}     
\def\la{{\triangleright}} 
\def\id{{\rm id}} 
\def\eps{\epsilon}     
\def\span{{\rm span}}     
\def\q2{{q^{-2}}}     
\def\bicross{{\blacktriangleright\!\!\!\triangleleft}}     
\def\note#1{{}}     
\def\eqn#1#2{\begin{equation}#2\label{#1}\end{equation}}     
\def\qbinom#1#2#3{\left(\begin{array}{c}#1\\#2\end{array}\right)\sb#3}     
\def\Z{{\Bbb Z}}     
\def\can{{\rm can}}     
\def\note#1{}     
\def\CA{{\frak A}}  
\def\tA{{\tilde{A}}}
\def\tC{{\tilde{C}}}
\def\tpsi{{\tilde{\psi}}}
\def\mea#1#2{{\mid\!\!\!-\!\!\!{#1\over#2}\!\!\!\!-}}

\noindent{Running title:} {\em Modules and coalgebra Galois extensions}
\vspace{.9in}     
\begin{center}     
{\large\bf ON MODULES ASSOCIATED TO COALGEBRA GALOIS EXTENSIONS}     
\vspace{12pt}\\     
Tomasz Brzezi\'nski\footnote{Lloyd's Tercentenary Fellow.     
On~leave from:      
Department of Theoretical Physics, University of \L\'od\'z,     
Pomorska 149/153, 90--236 \L\'od\'z, Poland.     
{\sc e-mail}: {\tt tb10@york.ac.uk}}
\vspace{6pt}\\     
\normalsize  Department of Mathematics \\ University
of York, Heslington, York YO1 5DD, U.K.
\end{center}     
\vspace{24pt}

\begin{abstract}     
For a given entwining structure $(A,C)_\psi$ 
involving an algebra $A$, a coalgebra
$C$, and an entwining map $\psi: C\otimes A\to A\otimes C$, a category
$\M_A^C(\psi)$ of right $(A,C)_\psi$-modules is defined and its structure
analysed. In particular, the notion of a measuring of $(A,C)_\psi$ to
$(\tA,\tC)_\tpsi$ is introduced, and certain functors between
$\M_A^C(\psi)$ and $\M_\tA^\tC(\tpsi)$ induced by such a
measuring are defined. It is shown that these functors  are 
inverse equivalences iff they are exact (or one of them faithfully
exact) and the measuring satisfies a certain Galois-type condition. Next, left modules $E$ and right
modules $\bar{E}$ associated to a $C$-Galois extension $A$ of $B$ are   
defined. These can be thought of as objects dual to fibre bundles with   
coalgebra $C$ in the place of a structure group, and a fibre
$V$. Cross-sections of such    
associated modules are defined as module maps $E\to B$ or $\bar{E}\to 
B$. It is shown    
that they can be identified with suitably equivariant maps from the   
fibre to $A$. Also, it is shown that a $C$-Galois extension is cleft   
if and only if $A=B\tens C$ as left $B$-modules and right   
$C$-comodules. The relationship between the modules $E$ and $\bar{E}$ is
studied in the case when $V$ is finite-dimensional and in the case
when the canonical entwining map is bijective. 
\end{abstract}
~\vspace{.5in}

     
\newpage    

\section*{\normalsize\sc\centering 1. Introduction}\vspace{-.6\baselineskip}  
\setcounter{section}{1}  

The notion of a Hopf-Galois extension arose from the works of Chase   
and Sweedler \cite{ChaSwe:hop} and Kreimer and Takeuchi    
\cite{KreTak:hop} (see \cite{Mon:hop} for a review). From the   
geometric point of view,    
a Hopf-Galois extension is a dualisation of the notion of a principal   
bundle and  thus it is a cornerstone of  the Hopf algebra or quantum  
group gauge    
theory. Such a gauge theory, in the sense of connections, gauge   
transformations, curvature etc. on Hopf-Galois   
extensions was proposed in \cite{BrzMa:gau} and later developed in   
\cite{Haj:str} \cite{BrzMa:dif} \cite{Brz:tra}. Also, the notion of a   
quantum fibre bundle as a module associated to the Hopf-Galois   
extension was introduced in \cite{BrzMa:gau}. This led to quantum   
group version of objects important in classical gauge theory such as   
sections of a vector bundle. Slightly different   
approaches to quantum group gauge theory, which take principal  
bundles as a framework of such a theory but do not use the Hopf-Galois   
extensions explicitly,   were also proposed in   
\cite{Dur:geo}, \cite{Pfl:fib}.   
     
Motivated by the structure of quantum homogeneous spaces, the notion of    
a $C$-Galois extension $A$ of an algebra $B$ was recently introduced  
\cite{BrzMa:coa} as an object dual to     
a principal bundle on such a space. This has been done   
by requiring that $A$ and $C$ admit an  {\em entwining structure}  
specified by a map $\psi:C\tens A\to A\tens C$ satisfying a set of  
(self-dual) conditions (cf. Definition~\ref{ent}).    
A gauge theory, in   
the above sense, on such a $C$-Galois extension was developed. In the   
present paper we derive the algebraic version of the classical
correspondence between the gauge transformations (vertical automorphisms)
and ad-covariant functions on a principal bundle (sections of an
associated adjoint bundle). The main objective of the present paper,
however, is to construct the algebraic counterpart of the notion of an
associated fibre bundle - a ``coalgebra fibre bundle". Our construction 
is motivated by a recent development of 
quantum and braided group Riemannian geometry in \cite{Maj:rie} and is 
a starting point for a more general coalgebra Riemannian geometry 
which is presented in \cite{BrzMa:geo}. The idea of the construction is
to associate a certain $B$-module to a
$C$-Galois extension of $B$ and a $C$-comodule. There are two possibilities of
associating such modules:  they  can be either left or
right $B$-modules, depending on whether there is a left or right
$C$-comodule involved (as opposed to the Hopf-Galois case, where the similar
construction leads to bimodules). We study both cases separately as well as the
relationship between them. In both cases we derive the algebraic
counterparts of the classical geometric equivalences between 
cross-sections and equivariant functions on a fibre bundle, and
between  cross-sections and trivialisations of a principal bundle,
and thus we generalise the Hopf-Galois considerations of \cite{Brz:tra}
to the $C$-Galois case.
It turns out that to perform this analysis it is useful to
consider the 
category $\M_A^C(\psi)$ of (right) $(A,C)_\psi$-modules. These are  
a natural generalisation of right $(A,H)$-Hopf modules.
We introduce the notion of a measuring of entwining structures, and
study when the functors between categories of entwined modules induced
by such a measuring are inverse equivalences, thus extending the results 
of \cite{CaeRai:ind} proven  for a
generalisation of Hopf modules known as Doi-Hopf modules \cite{Doi:uni}
\cite{Kop:var}.
\vspace{.3\baselineskip}  
   
\noindent {\sc Notation.} We work over a ground field    
$k$. All algebras are      
 associative and unital with the unit denoted by~1 (the unit map from
$k$ to the algebra is
denoted by $\eta$). We use the standard     
algebra and coalgebra notation, i.e., $\Delta$ is a coproduct, $\mu$ is a     
product, $\eps$ is a counit, etc. The identity map from the space $V$   
to itself      
is also denoted by $V$. The unadorned tensor product stands for the tensor     
product over $k$.     
For an algebra $A$ we denote by $\M_A$     
(resp. $_A\M$) the category of  right (resp. left) $A$-modules. 
For a right (resp. left)  
$A$-module $V$ the action is denoted by $\mu_V$ (resp. $_V\mu$) whenether it  
needs to be specified as a map, or by a dot between  
elements. Similarly, for a      
coalgebra $C$  we      
 denote by $\M^C$ (resp. $^C\M$) the category of  right  (resp. left)   
$C$-comodules.     
A right (resp. left)  coaction of $C$ on $V$ is denoted by $\Delta_V$  
(resp. ${}_V\Delta$). Also, by $_A\CM^C$ we     
denote the category of $(A,C)$-bimodules, {\em i.e.} left     
$A$-modules   and right      
$C$-comodules $V$  such that     
$\Delta_V\circ{}_V\mu = ({}_V\mu\tens C)\circ (A\tens\Delta_V)$,     
i.e., $\Delta_V$  is left  $A$-linear.  Similarly $^C\CM_A$ is the
category of right $A$-modules and left  
$C$-comodules with right $A$-linear coaction.  The position of a 
subscript (resp. superscript) of $\rm Hom$, $\rm End$ etc. indicates
left or right module (resp. comodule) structure, e.g. ${\rm Hom}_{-A}$
are morphsisms in $\M_A$. 
For coactions and coproducts we use Sweedler's notation with suppressed
summation sign: $\Delta_A (a) = a\sw 0\tens a\sw 1$ for $a\in A\in\M^C$,  
$\;_V\Delta(v) = v\sw{-1}\tens v\sw 0$ for $v\in V\in{}^C\M$ and     
$\hD(c)=c\1\ot c\2\,$ for $c\in C$. For an algebra $A$ and a coalgebra
$C$ we denote by $*$ the {\em convolution product} in ${\rm
Hom}(C,A)$, i.e. $f*g(c) = f(c\sw 1)g(c\sw 2)$ for any $f,g\in {\rm
Hom}(C,A)$ and $c\in C$. The convolution product makes ${\rm
Hom}(C,A)$ into an associative algebra with unit $\eta\circ \eps$. 
An element $f\in {\rm
Hom}(C,A)$ is said to be {\em convolution invertible} if it is
invertible with respect to $*$. 
     
\section*{\normalsize\sc\centering 2. $C$-Galois extensions and  
their automorphisms}\vspace{-.6\baselineskip}      
\setcounter{section}{2}     
  
First recall the definition of a coalgebra Galois extension from  
\cite{BrzHaj:coa}   
     
\begin{definition}  
Let $C$ be a coalgebra,     
$A$ an algebra and  a right $C$-comodule, and $B$ a subalgebra of $A$,  
$B :=\{b\in A\,|\;\forall a\in A \; \hD_A(ba)=ba\sw 0\otimes a\sw  
1\}$. We say that      
A is a {\em coalgebra Galois extension}      
(or $C$-Galois extension) of $B$ iff the canonical left $A$-module  
right $C$-comodule map      
\[     
can:=(\mu\ot C)\ci (A\ot\sb B \hD_A)\, :\; A\ot\sb B A\lra A\ot  
C\vspace*{-2.5mm}      
\]     
is bijective.  
Such a $C$-Galois extension is denoted by $A(B)^C$.    
\label{cge}  
\end{definition}     
  
Definition~\ref{cge}  
generalises the notion of a Hopf-Galois extension   
(see \cite{Mon:hop} for a review). The latter is a $C$-Galois
extension with $C=H$ being a Hopf algebra and $A$ a right $H$-comodule
algebra.   

An important role in the  
analysis of  coalgebra Galois extensions is   
played by the notion of an entwining structure \cite{BrzMa:coa} (closely   
connected with the theory of factorisation of algebras \cite{Maj:phy}).

\bde\label{ent}     
Let $C$ be a coalgebra, $A$ an algebra and let $\psi$ be a     
$k$-linear map $\psi: C\tens A\to A\tens C$ such that     
\begin{equation}\label{diag.A}     
\psi\circ(C\tens \mu) = (\mu\tens C)\circ (A\tens\psi)\circ(\psi\tens A),     
\qquad \psi\circ (C\tens\eta) = \eta\tens C,     
\end{equation}     
\begin{equation}\label{diag.B}     
(A\tens\Delta)\circ\psi = (\psi\tens     
C)\circ(C\tens\psi)\circ(\Delta\tens A), \qquad (A\tens \eps)\circ\psi =     
\eps\tens A,     
\end{equation}     
The triple     
$(A,C,\psi)$ is called an {\em entwining structure} and is
denoted by $(A,C)_\psi$.  The map $\psi$ is called an {\em entwining
map}. A morphism between entwining structures
$(A,C)_\psi$ and $(\tA,\tC)_\tpsi$ is a pair $(f,g)$, where $f:A\to
\tA$ is a unital algebra map and $g:C\to \tC$ is a countial coalgebra
map such that $(f\otimes g)\circ\psi = \tpsi\circ(g\otimes f)$.
\ede

Given an entwining structure $(A,C)_\psi$ we 
use the notation $\psi(c\tens a) = a_\alpha\tens c^\alpha $     
(summation over a Greek index is understood), for all
$a\in A$, $c\in C$.

For an entwining structure $(A,C)_\psi$, $\M_A^C(\psi)$ is 
the category of {\em right $(A,C)_\psi$-modules}. The objects of  
$\M_A^C(\psi)$ are right $A$-modules and right $C$-comodules $M$ such  
that 
 \begin{equation}\label{more.coa}  
\Delta_A(m\cdot a) = m\sw 0\cdot\psi(m\sw 1\tens a) := m\sw 0\cdot  
a_\alpha \tens m\sw 1^\alpha, \qquad\qquad \forall m\in M,a\in A  
\end{equation}  
Morphisms in $\M_A^C(\psi)$ are right $A$-module right  
$C$-comodule maps. 

It is shown in \cite{BrzHaj:coa} that if $A(B)^C$ is a $C$-Galois 
extension, then $\psi: C\tens A\to A\tens C$, $\psi = can     
\circ(A\otimes_B\mu)\circ(\tau \tens A)$ is a unique entwining map
such that $A$ is an object in $\M_A^C(\psi)$.
Here $\tau : C\to A\otimes_B A$ is the {\em translation map},     
i.e. $\tau(c) = can^{-1}(1\otimes c)$. This $(A,C)_\psi$ is called the {\em     
canonical entwining structure} associated to 
$A(B)^C$.    
 
A $C$-Galois extension  $A(B)^C$ is called
a {\em cleft extension} iff there exists a convolution invertible, right
$C$-comodule map $\Phi: C\to A$. Such a $\Phi$ is called a {\em cleaving map}.
An object dual to a trivial principal bundle is an example of a 
cleft extension. The following proposition gives
equivalent descriptions of cleft extensions, generalising
\cite[Proposition~2.9]{BrzMa:coa}, where more structure on $C$ and a
different condition for $\Phi$ were assumed.

\bpr\label{ex.cleft}
 Let $C$ be a coalgebra, $A$ be a right $C$-comodule and let $B$
be as in Definition~\ref{cge}. If there exists a convolution
invertible,  right $C$-comodule map   
$\Phi:C\to A$ 
then the following are equivalent:

(1) $A$ is a $C$-Galois extension of $B$.

(2) There is an entwining structure $(A,C)_\psi$ such that $A\in
\M_A^C(\psi)$ via $\Delta_A$ and $\mu$.

(3) For every $a\in A$, $a\sw 0\Phi^{-1}(a\sw 1)\in B$

If any of the above conditions hold, then  $A\cong B\tens C$   
in ${}_B\CM^C$. 
\epr
\proof (1) $ \Rightarrow$ (2) follows from the result of
\cite{BrzHaj:coa}, cited above. 

(2) $\Rightarrow$ (3) Since $A\in \M_A^C(\psi)$, the coaction can be
written as $\Delta_A(a) = 1\sw 0\psi(1\sw 1\otimes a)$, for all $a\in
A$. The right $C$-colinearity of $\Phi$ together with the equality 
$1\sw 0\eps(c)\tens 1\sw 1 = 1\sw 0\psi(1\sw 1\tens\Phi(c\sw  
1)\Phi^{-1}(c\sw 2))$ and (\ref{diag.A}) imply that
\begin{equation}
\psi(c\sw 1\tens \Phi^{-1}(c\sw 2)) = \Phi^{-1}(c)\Delta_A(1).   
\label{phi-1}
\end{equation}
Using this last equality and the fact that $A\in \M_A^C(\psi)$ one
finds that  (3) holds.

(3) $\Rightarrow$ (1) One easily verifies that the map 
$A\otimes C\to A\otimes_B A$,
$a\otimes c\mapsto a\Phi^{-1}(c\sw 1)\otimes_B\Phi(c\sw 2)$ is the
inverse of $can$.

Finally, a ${}_B\CM^C$ isomorphism $A\stackrel{\sim}{\to} B\otimes C$ is $a\mapsto a\sw
0\Phi^{-1}(a\sw 1)\otimes a\sw 2$, and its inverse is $b\otimes c\mapsto
b\Phi(c)$.
\endproof    

For $A(B)^C$, ${\rm Aut
}(A(B)^C)$ denotes the  group of left $B$-module, right
$C$-comodule automorphisms of $A$,
with the product ${\cal FG} = {\cal G}\circ {\cal F}$, for all ${\cal
F},{\cal G}\in {\rm Aut
}(A(B)^C)$. We  
now give a description of ${\rm Aut
}(A(B)^C)$ which reflects the classical correspondence between
gauge transformations of a principal bundle and 
ad-equivariant functions on it 
(cf. \cite[7.1.6]{Hus:fib})\footnote{In the Hopf-Galois case, 
(\ref{ad-cov}) means that $f\in{\rm Hom}^{-C}(C,A)$, where $C$ is
in $\M^C$ via the right adjoint coaction. Also, a condition similar to 
(\ref{ad-cov}) characterises connection one-forms on $A(B)^C$.}.
\begin{theorem}\label{C-aut}  
The group ${\rm Aut
}(A(B)^C)$ is isomorphic to the group
$\CC(A)$ of convolution invertible maps $f:C\to A$ such that   
\begin{equation}\label{ad-cov}  
\psi\circ(C\tens f)\circ\Delta = (f\tens C)\circ\Delta,  
\end{equation}  
where $\psi$ is the canonical entwining map. The product in $\CC(A)$
is the convolution product.  
\end{theorem}  
\proof  We use the $\pi$-method of \cite{DoiTak:miy}. Applying the
functor ${\rm Hom}_{A-}(-,A)$ to $can : A\otimes_B A\stackrel{\sim}{\to}
 A\otimes C$
one obtains the isomorphism $\pi : {\rm Hom} (C,A)\stackrel\sim\to {\rm
Hom}_{B-}(A,A)$. Explicitly, $\pi(f) = \mu\circ(A\otimes
f)\circ\Delta_A$ and $\pi^{-1}({\cal F}) = \mu\circ(A\otimes_B{\cal F})\circ\tau$, 
for all $f\in {\rm Hom} (C,A)$, ${\cal F}\in {\rm
Hom}_{B-}(A,A)$.  Since $A\in \M_A^C(\psi)$ we have $\Delta_A(\pi(f)(a)) =  \Delta_A(a\sw 0
f(a\sw 1))   
= a\sw 0\psi(a\sw 1\tens f(a\sw 2))$. It means that $\pi(f)$ is right $C$-colinear if and
only if for all $a\in A$, 
$
a\sw 0\psi(a\sw 1\tens f(a\sw 2)) = a\sw 0f(a\sw 1)\tens a\sw 2 .
$
Using the definition of $\tau$ one easily sees that this is equivalent to
(\ref{ad-cov}). Next take any $f,g\in 
\CC(A)$. Since $\pi(f)$ is a right $C$-comodule map, we find for all
$a\in A$
$$
(\pi(g)\circ\pi(f))(a) = \pi(f)(a\sw 0)g(a\sw 1) = a\sw 0f(a\sw 1)g(a\sw
2) = \pi(f*g)(a).
$$
Finally, if $f$ satisfies  (\ref{ad-cov}) then, 
using (\ref{diag.A}), for all  
$c\in C$ one has   
$$
c\sw 1\tens 1\tens c\sw 2 = c\sw 1\tens\psi(c\sw 2\tens f(c\sw  
3)f^{-1}(c\sw 4)) 
= c\sw 1 \tens f(c\sw 2)\psi(c\sw 3\tens f^{-1}(c\sw  
4)).  
$$
Applying $f^{-1}\tens A\tens C$ to the above equality and multiplying  
the first two factors one finds that $f^{-1}$ satisfies (\ref{ad-cov}). 
Also $\eta\circ \eps$ 
satisfies  (\ref{ad-cov}). Therefore $\CC(A)$ is a group with respect
to the convolution product as claimed, and 
${\rm Aut }(A(B)^C) \cong \CC(A)$ as groups. \endproof  
 
Finally we notice that if  $A(B)^C$ is a cleft $C$-Galois extension,
then since $A\cong B\tens C$   
in ${}_B\CM^C$,  there is an algebra isomorphism ${\rm
End}_{B-}^{-C}(A)
\cong {\rm Hom}(C,B)^{op}$ (cf. \cite[Lemma on p.91]{Mon:hop}). This
implies that, in this case, the  group ${\rm Aut }(A(B)^C)$ and
therefore also $\CC(A)$ are 
isomorphic to the group of convolution   
invertible maps $C\to B$. This reflects the classical description of
local gauge transformations (cf. \cite[7.1.7]{Hus:fib}).
 
\section*{\normalsize\sc\centering 3. The structure of  
$(A,C)_\psi$-modules}\vspace{-.6\baselineskip}    
\setcounter{section}{3}     
\setcounter{definition}{0}  
In this section we analyse the structure of the category  
$\M_A^C(\psi)$ of $(A,C)_\psi$-modules, i.e. right $A$-modules and  
right $C$-comodules characterised by (\ref{more.coa}). This category  
can  
be viewed as a generalisation of the categories well-studied
in the Hopf algebra theory.

\bex\label{ex.cat}{\rm 

(1) Let $C=H$ be a Hopf algebra, $A$ be a right $H$-comodule algebra
and let $\psi: H\otimes A\to A\otimes H$ be defined by $\psi: h\otimes
a\mapsto a\sw 0\otimes ha\sw 1$. Then $\M_A^H(\psi)$ is the category
of right $(A,H)$-Hopf modules introduced in \cite{Doi:rel}.

(2) Let $A=C=H$ be a Hopf algebra and let the entwining map
$\psi:H\otimes H\to H\otimes H$ be given by $\psi: g\otimes h \mapsto
h\sw 2 \otimes (Sh\sw 1)gh\sw 3$, where $S$ is the antipode in $H$. 
Then $\M_H^H(\psi)$ is the category
of right-right Yetter-Drinfeld modules introduced in
\cite{Yet:rep}, \cite{RadTow:yet}.

(3) Examples (1) and (2) are special cases of the following
construction. Let $H$ be a Hopf algebra, $A$ be a right $H$-comodule
algebra and $C$ a right $H$-module coalgebra. Then $(A,C)_\psi$ is an
entwining structure with $\psi: c\otimes a \mapsto a\sw 0\otimes c\cdot 
a\sw 1$ and $\M_A^C(\psi)$ is the category of unifying Hopf modules (or
Doi-Hopf modules) introduced in \cite{Doi:uni} \cite{Kop:var}. A special
case of this category with $C = H/I$ a quotient
coalgebra and a quotient right $H$-module was considered in
\cite{Sch:pri} (in particular, $\M_H^C(\psi)$ for the canonical
entwining structure $(H,C)_\psi$ associated to a $C$-Galois extension 
$H(B)^C$, where $B$ is a quantum homogeneous space of $H$ 
(cf. \cite[Example~2.5]{BrzMa:coa} or 
\cite[Lemma~1.3.]{Sch:nor}), is of this type).}
\eex

\bde\label{measure}
Let $(A,C)_\psi$ and $(\tA,\tC)_\tpsi$ be entwining structures, and
let $\alpha :\tC\otimes\tA\to A$, $\gamma:\tC\to A\otimes C$ be linear
maps such that
\begin{equation}
\mu\circ(\alpha\otimes\alpha)\circ (\tC\otimes \tpsi\otimes
\tA)\circ(\tilde{\Delta}\otimes \tA\otimes \tA) = \alpha\circ(\tC\otimes
\tilde{\mu}),
\qquad \alpha\circ(\tC\otimes \tilde{\eta}) = \eta\circ\tilde{\eps},
\label{mea1}
\end{equation}
\begin{equation}
(\mu\otimes C\otimes C)\circ (A\otimes\psi\otimes C)\circ(\gamma\otimes\gamma)
\circ\tilde{\Delta} = (A\otimes\Delta)\circ\gamma, \qquad
(A\otimes\eps)\circ\gamma = \eta\circ\tilde{\eps},
\label{mea2}
\end{equation}
\begin{equation}
(\mu\otimes C)\circ(\alpha\otimes\gamma)\circ(\tC\otimes\tpsi)
\circ(\tilde{\Delta}\otimes \tA) = (\mu\otimes C)\circ(A\otimes
\psi)\circ(\gamma\otimes\alpha)\circ(\tilde{\Delta}\otimes \tA),
\label{mea3}
\end{equation}
where non-tilded (tilded) structure maps correspond to $A,C$
($\tA,\tC$). The pair
$(\alpha,\gamma)$ is said to {\em measure} $(\tA,\tC)_\tpsi$ to
$(A,C)_\psi$. Such a measuring is denoted by $(\tA,\tC)_\tpsi
\mea\alpha\gamma (A,C)_\psi$.
\ede
The terminology of Definition~\ref{measure} is motivated by the fact
that if one chooses $\psi$ and $\tpsi$ to be the twists, and $C=k$, then
the pair $(\alpha,\eta\circ\tilde{\eps})$ measures $(\tA,\tC)_\tpsi$ to
$(A,k)_\psi$ iff $(\alpha, \tC)$ measures $\tA$ to $A$ in the sense of 
\cite[p.~138]{swe}. If $(f,g)$ is a morphism from  $(\tA,\tC)_\tpsi$
to $(A,C)_\psi$ 
then $(\tilde{\eps}\otimes f, \eta\otimes g)$ measures 
$(\tA,\tC)_\tpsi$ to
$(A,C)_\psi$.

\bpr\label{functors}
Let $(\alpha,\gamma)$ measure  $(\tA,\tC)_\tpsi$ to
$(A,C)_\psi$. Then:

(1) For all $M\in \M_A$, $M\otimes \tC$ is an
$(\tA,\tC)_\tpsi$-module via $\Delta_{M\otimes\tC} =
M\otimes\tilde{\Delta}$ and $\mu_{M\otimes\tC} =
(\mu_M\otimes\tC)\circ(M\otimes\alpha\otimes
\tC)\circ(M\otimes\tC\otimes\tpsi)\circ(M\otimes\tilde{\Delta}
\otimes\tA)$.

(2) For all $\tilde{M}\in\M^\tC$, $\tilde{M}\otimes
A$ is an $(A,C)_\psi$-module via
$\mu_{\tilde{M}\otimes A} =
\tilde{M}\otimes\mu $ and $\Delta_{\tilde{M}\otimes A}
= (\tilde{M}\otimes\mu\otimes
C)\circ(\tilde{M}\otimes
A\otimes\psi)\circ(\tilde{M}\otimes\gamma\otimes
A)\circ(\Delta_{\tilde{M}}\otimes A)$.
\epr
\proof (1) We first show that  $\mu_{M\otimes\tC}$ (later denoted by a
dot) 
is an action of $\tA$ on $M\otimes \tC$. Explicitly, 
for any $m\in M$, $\tilde{c}\in\tC$ and $\tilde{a}\in \tA$ this map is
$(m\otimes \tilde{c})\cdot\tilde{a} = m\cdot\alpha(\tilde{c}\sw
1\otimes\tilde{a}_\beta)\otimes \tilde{c}\sw 2^\beta$. By the second
of equations (\ref{mea1}) and (\ref{diag.A}) we have that 
$(m\otimes \tilde{c})\cdot 1 = m\otimes
\tilde{c}$. Furthermore, for any
$\tilde{a}'\in \tA$,
\begin{eqnarray*}
((m\otimes \tilde{c})\cdot\tilde{a})\cdot\tilde{a}' & = & 
m\cdot\alpha(\tilde{c}\sw
1\otimes\tilde{a}_\beta)\alpha(\tilde{c}\sw 2^\beta\sw 1 \otimes
\tilde{a}'_\delta) \otimes \tilde{c}\sw 2^\beta\sw 2^\delta \\
& = & m\cdot\alpha(\tilde{c}\sw
1\otimes\tilde{a}_{\beta\lambda})\alpha(\tilde{c}\sw 2^\lambda \otimes
\tilde{a}'_\delta) \otimes \tilde{c}\sw 3^{\beta\delta} \qquad \qquad\mbox{\rm
(by (\ref{diag.B}))} \\
& = & m\cdot\alpha(\tilde{c}\sw
1\otimes\tilde{a}_\beta\tilde{a}'_\delta)\otimes \tilde{c}\sw
2^{\beta\delta} \qquad\qquad\qquad\qquad\quad \mbox{\rm (by (\ref{mea1}))}\\
& = & (m\otimes \tilde{c})\cdot(\tilde{a}\tilde{a}')
\qquad\qquad\qquad\qquad\qquad\qquad\qquad  \mbox{\rm
(by (\ref{diag.A}))}.
\end{eqnarray*}
Clearly, $\Delta_{M\otimes \tC}$ is a right coaction of $\tC$ on
$M\otimes \tC$. For any $m\in M$, $\tilde{c}\in\tC$, $\tilde{a}\in \tA$, 
\begin{eqnarray*}
\Delta_{M\otimes \tC}((m\otimes\tilde{c})\cdot \tilde{a}) & = &
m\cdot\alpha(\tilde{c}\sw 1\otimes \tilde{a}_\beta)\otimes \tilde{c}\sw
2\sp\beta\sw 1\otimes \tilde{c}\sw 2\sp\beta\sw 2\\
& = & m\cdot\alpha(\tilde{c}\sw 1\otimes \tilde{a}_{\beta\delta})
\otimes \tilde{c}\sw
2\sp\delta\otimes \tilde{c}\sw 3\sp\beta \qquad \mbox{\rm
(by (\ref{diag.B}))} \\
& = & (m\otimes \tilde{c}\sw 1)\cdot \tilde{a}_\beta\otimes \tilde{c}\sw
2^\beta.
\end{eqnarray*}
This proves that $M\otimes\tC$ is an object in $\M_\tA^\tC(\tpsi)$.

(2) Dual to (1).
\endproof

\bco\label{MC1} Let $(A,C)_\psi$ be an entwining structure. Then:  
  
(1) For any right $A$-module $M$, $M\otimes C$ is an  
$(A,C)_\psi$-module with the action $m\otimes c \otimes a \mapsto  
m\cdot\psi(c\otimes a)$ and the coaction $\Delta_{M\otimes C} =  
M\otimes \Delta$.  
  
(2) For any right $C$-comodule $V$, $V\otimes A$ is an  
$(A,C)_\psi$-module with the action $V\otimes \mu$ and the coaction  
$v\otimes a \mapsto v\sw 0\otimes \psi(v\sw 1\otimes a)$.  
\eco  
\proof To prove (1) take 
$(A,k)_\sigma$, where $\sigma :k\otimes A\to A\otimes k$ is a twist 
(canonically equivalent to
the map $A$) and notice that $
(\eps\otimes A,\eta\circ\eps)$, measures $(A,C)_\psi$ to $(A,k)_\sigma$.
Then  Proposition~\ref{functors}(1) yields the assertion. 
Statement (2) is
dual to (1), and can be deduced from Proposition~\ref{functors}(2) by
taking $(\tA,\tC)_\tpsi = (k,C)_{\sigma}$ and 
$(\alpha,\gamma) = (\eta\circ\eps, \eta\otimes C)$. \endproof

\bpr\label{maps} Let $(\alpha,\gamma)$ measure  $(\tA,\tC)_\tpsi$ to 
$(A,C)_\psi$. Then:

(1) For all $M\in \M_A^C(\psi)$ the map $\hat{\ell}^M:M\otimes \tC\to
M\otimes C\otimes \tC$, 
$$
\hat{\ell}^M := \Delta_M\otimes \tC - (\mu_M\otimes C\otimes
\tC)\circ(M\otimes\gamma\otimes \tC)\circ(M\otimes\tilde{\Delta})
$$ is a
morphism in $\M_\tA^\tC(\tpsi)$, where $M\otimes \tC$ and $M\otimes
C\otimes \tC$ are viewed in $\M_\tA^\tC(\tpsi)$ as in
Proposition~\ref{functors}(1) with $M\otimes C\in \M_A$ as in
 Corollary~\ref{MC1}(1).

(2) For all $\tilde{M}\in \M_\tA^\tC(\tpsi)$, the map
$\hat\ell_{\tilde{M}} : \tilde{M}\otimes \tA\otimes A\to \tilde{M}\otimes
A$, 
$$
\hat\ell_{\tilde{M}} = \mu_{\tilde{M}}\tens A - (\tilde{M}\tens
\mu)\circ (\tilde{M}\tens \alpha\tens A)\circ(\Delta_{\tilde{M}}\tens
\tA\tens A)
$$
is a morphism in $\M_A^C(\psi)$, where $\tilde{M}\otimes A$ and $\tilde{M}\otimes
\tA\otimes A$ are viewed in $\M_A^C(\psi)$ as in
Proposition~\ref{functors}(2) with $\tilde{M}\otimes \tA\in \M^\tC$ as
in Corollary~\ref{MC1}(2).
\epr
\proof We introduce the Sweedler-like notation $\gamma(\tilde{c}) =
\tilde{c}\sua{1}\otimes \tilde{c}\sua{2}$ (summation understood). With
this notation the map $\hat{\ell}^M$ explicitly reads for all $m\in M$,
$\tilde{c}\in \tC$, $\hat{\ell}^M(m\otimes \tilde{c}) = m\sw 0\otimes m\sw
1\otimes \tilde{c} - m\cdot \tilde{c}\sw 1\sua 1\otimes \tilde{c}\sw
1\sua 2\otimes \tilde{c}\sw 2$. Clearly, $\hat{\ell}^M$ is a right 
$\tC$-comodule map. Next we have
\begin{eqnarray*}
\hat{\ell}^M((m\otimes\tilde{c})\cdot \tilde{a}) & = &
\hat{\ell}^M(m\cdot\alpha(\tilde{c}\sw 1\otimes \tilde{a}_\beta)\otimes
\tilde{c}\sw 2^\beta) \\
& = & m\sw 0\cdot\alpha(\tilde{c}\sw 1\otimes \tilde{a}_\beta)_\delta
\otimes m\sw 1^\delta \otimes \tilde{c}\sw 2^\beta \qquad\qquad\qquad 
\mbox{\rm ($M\in \M_A^C(\psi)$)} \\
&&- m\cdot
\alpha(\tilde{c}\sw 1\otimes \tilde{a}_\beta)\tilde{c}\sw 2^\beta\sw 1\sua
1\otimes   \tilde{c}\sw 2^\beta\sw 1\sua 2\otimes \tilde{c}\sw 2^\beta\sw
2 \\
& = & m\sw 0\cdot\alpha(\tilde{c}\sw 1\otimes \tilde{a}_\beta)_\delta
\otimes m\sw 1^\delta \otimes \tilde{c}\sw 2^\beta\\
&& - m\cdot
\alpha(\tilde{c}\sw 1\otimes \tilde{a}_{\beta\delta})\tilde{c}\sw
2^\delta\sua
1\otimes   \tilde{c}\sw 2^\delta\sua 2\otimes \tilde{c}\sw 3^\beta 
\qquad\qquad\quad \mbox{\rm
(by (\ref{diag.B}))}\\
& = & m\sw 0\cdot\alpha(\tilde{c}\sw 1\otimes \tilde{a}_\beta)_\delta
\otimes m\sw 1^\delta \otimes \tilde{c}\sw 2^\beta\\
&& - m\cdot
\tilde{c}\sw 1\sua 1\alpha(\tilde{c}\sw 2\otimes \tilde{a}_\beta)_\delta
\otimes \tilde{c}\sw 1\sua 2^\delta\otimes \tilde{c}\sw 3^\beta \qquad
\qquad \quad \mbox{\rm
(by (\ref{mea3}))}\\
& = & (m\sw 0\otimes m\sw 1)\cdot\alpha(\tilde{c}\sw 1\otimes 
\tilde{a}_\beta)\otimes \tilde{c}\sw 2^\beta\\
&& - (m\cdot
\tilde{c}\sw 1\sua 1\otimes \tilde{c}\sw 1\sua 2)\cdot
\alpha(\tilde{c}\sw 2\otimes \tilde{a}_\beta)
\otimes \tilde{c}\sw 3^\beta\\
& = & (m\sw 0\otimes m\sw 1\otimes \tilde{c})\cdot 
\tilde{a} - (m\cdot
\tilde{c}\sw 1\sua 1\otimes \tilde{c}\sw 1\sua 2\otimes \tilde{c}\sw 2)
\cdot \tilde{a}\\
& = & \hat{\ell}^M(m\otimes \tilde{c})\cdot \tilde{a}.
\end{eqnarray*}
To derive the fifth and the sixth equations we used definitions of actions
of $A$ on $M\otimes C$ in Corollary~\ref{MC1}(1) and of $\tA$ on $M\otimes
C\otimes \tC$ in Proposition~\ref{functors}(1) combined with
Corollary~\ref{MC1}(1). This completes the proof that $\hat{\ell}^M$ is a
morphism in $\M_\tA^\tC(\tpsi)$.

(2) Dual to (1).
\endproof

Given a measuring $(\tA,\tC)_\tpsi\mea\alpha\gamma
(A,C)_\psi$,  for all $M\in \M_A^C(\psi)$ define 
$M\widehat{\square}_C\tC\subset M\otimes \tC$ via the exact sequence
$$
\begin{CD}
0 @>>> M\widehat{\square}_C\tC @>>> M\otimes \tC @>{\hat{\ell}^M}>>
 M\otimes C\otimes\tC.
\end{CD}
$$
Since the above sequence is a sequence in $\M_\tA^\tC(\tpsi)$,
$M\widehat{\square}_C\tC$ is an $(\tA,\tC)_\tpsi$-module via the
restriction of the structure maps in Proposition~\ref{functors}(1). Thus we
obtain a functor $-\widehat{\square}_C\tC:\M_A^C(\psi)\to \M_\tA^\tC(\tpsi)$.
Dually, for all $\tilde{M}\in \M^\tC_\tA(\tpsi)$ define
$\tilde{M}\hat{\otimes}_\tA A $ by the exact sequence
$$
\begin{CD}
\tilde{M}\tens \tA\tens A @>{\hat\ell_{\tilde{M}}}>> \tilde{M}\tens A @>{\hat{\pi}_{\tilde{M}}}>>
\tilde{M}\hat{\tens}_\tA A @>>> 0
\end{CD}
$$
(if there is no need to specify the module $\tilde{M}$
 we will write $\hat\pi$ for $\hat{\pi}_{\tilde{M}}$). 
$\tilde{M}\hat{\otimes}_\tA A $ is an $(A,C)_\psi$-module with the
structure maps obtained from the structure maps in Proposition~\ref{functors}(2), 
by projecting through $\hat{\pi}_M$. Thus we have the functor 
$-\hat{\tens}_\tA
A: \M_\tA^\tC(\tpsi)\to M_A^C(\psi)$.
 
\bpr\label{prop.adjoint.n}
Given a measuring $(\tA,\tC)_\tpsi\mea{\alpha}{\gamma} (A,C)_\psi$,
the functor $-\widehat{\square}_C\tC:\M_A^C(\psi)\to \M_\tA^\tC(\tpsi)$
is the right adjoint of the functor $-\hat{\tens}_\tA A$.
\epr
\proof We claim that for all $M\in \M_A^C(\psi)$, 
$\tilde{M}\in\M_\tA^\tC(\tpsi)$ there
is a natural isomorphism $\zeta_{\tilde{M}, M}:{\rm
Hom}_{-A}^{-C}(\tilde{M}\hat{\tens}_\tA A, M)\to {\rm
Hom}_{-\tA}^{-\tC}(\tilde{M}, M\widehat{\square}_C\tC)$,
$\zeta_{\tilde{M},M}: f \mapsto (f\circ\hat\pi{\tens}
\tC)\circ(\tilde{M}\tens \eta\tens\tC)\circ\Delta_{\tilde{M}}$. 
Explicitly for  all $\tilde{m}\in \tilde{M}$, 
$\zeta_{\tilde{M},M}(f)(\tilde{m})= 
f(\hat\pi(\tilde{m}\sw 0\tens 1))\tens \tilde{m}\sw 1$.
 The output of $\zeta_{\tilde{M},M}(f)$ is in
$M\widehat{\square}_C\tC$, since
\begin{eqnarray*}
\hat{\ell}^M(\zeta_{\tilde{M},M}(f)(\tilde{m}))\!\!\!\!\! &= & \!\!\!\!\!
\Delta_M(f(\hat\pi(\tilde{m}\sw 0\tens 1)))\tens \tilde{m}\sw 1 
- f(\hat\pi(\tilde{m}\sw 0\otimes 1))\cdot \tilde{m}\sw 1\sua 1\otimes
\tilde{m}\sw 1\sua 2\otimes \tilde{m}\sw 2\\
& = &f(\hat\pi(\tilde{m}\sw 0\otimes 1)\sw 0)\otimes \hat\pi(\tilde{m}\sw
0\otimes 1)\sw 1\otimes \tilde{m}\sw 1 \\
&&- f(\hat\pi(\tilde{m}\sw 0\otimes 1)\cdot \tilde{m}\sw 1\sua 1)\otimes
\tilde{m}\sw 1\sua 2\otimes \tilde{m}\sw 2\\
& = & f(\hat\pi(\tilde{m}\sw 0\otimes \tilde{m}\sw 1\sua 1))\otimes
\tilde{m}\sw 1\sua 2\otimes \tilde{m}\sw 2 \\
&&- f(\hat\pi(\tilde{m}\sw 0\otimes \tilde{m}\sw 1\sua 1))\otimes
\tilde{m}\sw 1\sua 2\otimes \tilde{m}\sw 2\\
& = & 0,
\end{eqnarray*}
where we used the definition of $\hat{\ell}^M$ to derive the first
equality, then the fact that $f$ is a morphism in $\M_A^C(\psi)$ to
obtain the second one. The third equality was obtained by using the
explicit form of the coaction of $C$ on $\tilde{M}\hat{\otimes}_\tA A $ and the fact
that $\hat\pi$ is a right $A$-module map.

It is clear that $\zeta_{\tilde{M},M}(f)$ is a right $\tC$-comodule map,
it is also  right $\tA$-linear since
\begin{eqnarray*}
\zeta_{\tilde{M},M}(f)(\tilde{m})\cdot \tilde{a} & = & 
f(\hat\pi(\tilde{m}\sw 0\tens
1))\cdot \alpha(\tilde{m}\sw 1\otimes \tilde{a}_\beta)\tens
\tilde{m}\sw 2^\beta \\
& = & f(\hat\pi(\tilde{m}\sw 0\tens
\alpha(\tilde{m}\sw 1\otimes \tilde{a}_\beta)))\tens
\tilde{m}\sw 2^\beta  \qquad\qquad \mbox{\rm ($f$ is right $A$-linear)}\\
& = & f(\hat\pi(\tilde{m}\sw 0\cdot \tilde{a}_\beta\tens 1))
\tens
\tilde{m}\sw 1^\beta  \qquad\qquad\qquad\; \mbox{\rm (by definition of
$\tilde{M}\hat{\tens}_\tA A$)}\\
& = & \zeta_{\tilde{M},M}(f)(\tilde{m}\cdot \tilde{a})
\qquad\qquad\qquad\qquad\qquad\;\; \mbox{\rm ($\tilde{M}\in 
\M_\tA^\tC(\tpsi)$)}
\end{eqnarray*}
It is an easy exercise to verify that $\zeta_{\tilde{M},M}$ is natural
in $\tilde{M}$ and $M$ and that its inverse is $\zeta^{-1}_{\tilde{M},
M} (g)\circ \hat\pi =  \mu_M\circ(M\tens\tilde{\eps}\tens A)\circ(g\tens
A)$. 
\endproof

\bco\label{adjoint.cor} 
Let $(A,C)_\psi$ be an entwining structure. Then:  
  
(1) The functor $-\otimes C :\M_A\to \M_A^C(\psi)$ is the
right adjoint of the forgetful functor $\M_A^C(\psi)\to \M_A$.
  
(2) The functor $-\otimes A :\M^C\to \M_A^C(\psi)$ is the
left adjoint of the forgetful functor $\M_A^C(\psi)\to \M^C$.
\eco  
\proof To prove (1) take the measuring in the proof of Corollary~\ref{MC1}(1).
Then $\M_A^k(\sigma) = \M_A$, and for all $M\in \M_A$, $M\widehat{\square}_k C =
M\otimes C$, while for all $N\in\M_A^C(\psi)$, $N\hat{\otimes}_AA=N$. Now
Proposition~\ref{prop.adjoint.n} yields the assertion. 
Statement (2) is
dual to (1), and can be deduced from Proposition~\ref{prop.adjoint.n} by
taking $(\tA,\tC)_\tpsi = (k,C)_{\sigma}$ and 
$(\alpha,\gamma) = (\eta\circ\eps, \eta\otimes C)$. \endproof

From the proof of Proposition~\ref{prop.adjoint.n} it is clear that the
adjunction morphisms are 
$\Psi_{\tilde{M}} = (\hat\pi_{\tilde{M}}\otimes
\tC)\circ(\tilde{M}\otimes\eta\otimes \tC)\circ\Delta_{\tilde{M}}:\tilde{M} \to
(\tilde{M}\hat{\tens}_\tA
A)\widehat{\square}_C\tC$, and 
$\Phi_{M}: (M\widehat{\square}_C\tC)\hat{\tens}_\tA A\to M$ determined by
$\Phi_{M}\circ\hat\pi_{M\widehat{\square}_C\tC} = \mu_M\circ (M\otimes
\tilde{\eps}\tens A)$,
for all $M\in \M_A^C(\psi)$ and $\tilde{M}\in \M_\tA^\tC(\tpsi)$.

\begin{definition}{\rm (cf. Definition~1.4 in \cite{CaeRai:ind})}
$(\tA,\tC)_\tpsi \mea\alpha\gamma (A,C)_{\psi}$
is said to be a
{\em Galois measuring} if the adjunctions $\Psi_{\tC\otimes \tA}$,
$\Phi_{A\otimes C}$ are bijective ($C\otimes
A\in \M_A^C(\psi)$ and $\tilde{A}\otimes \tilde{C}\in
\M_{\tilde{A}}^{\tilde{C}}(\tilde{\psi})$ as in
Corollary~\ref{MC1}). 
\label{def.mea.g}
\end{definition}
\bex\label{just}
{\rm Assume that $A$ is an $(A,C)_\psi$-module and let $B := \{b\in A\; |\;
\forall a\in A, \;\Delta_A(ba) = b\Delta_A(a)\}$. 
Take the trivial entwining structure $(B,k)_\sigma$, so that
$\M_B^k(\sigma) = \M_B$. Then $(\iota_B, \Delta_A\circ\eta)$, 
where $\iota_B:B\hookrightarrow A$ is
the canonical inclusion, measures $(B,k)_\sigma$ to $(A,C)_\psi$. 
With this measuring, for all
$V\in \M_B$,
$V\hat{\tens}_B A = V\otimes_B A$, while for all $M\in \M_A^C(\psi)$,
$M\widehat{\square}_Ck = M_0 := \{m\in M \; | \; \Delta_M(m) = m\cdot 1\sw
0\otimes 1\sw 1\}$. Notice that  $  
M_0=\{m\in M\,|\;\forall  
a\in A, \; \hD_M(m\cdot a)=m\cdot a\sw 0\otimes a\sw   
1\}$. Indeed,    
if $m\in M_0$ then 
$$  
\Delta_M(m\cdot a)=m\sw 0\cdot \psi(m\sw 1 \tens a) = m\cdot (1\sw 0  
\psi(1\sw  1 \otimes a)) = m\cdot a\sw 0\otimes a\sw 1,  
$$  
since $\Delta_A(a) = 1\sw 0\psi(1\sw 1\otimes a)$. In
particular $B=A_0$, so that 
$\Psi_{B\otimes k}=B$. On
the other hand $A\cong (A\otimes C)_0$ via $a\mapsto a1\sw 0\otimes 1\sw
1$, $\sum_ia^i\otimes c^i \mapsto \sum_i \eps(c^i)a^i$. Taking this
isomorphism into account we have $\Phi_{A\otimes C} = can$, and we
conclude that $(B,k)_\sigma\mea{\iota_B}{\Delta_A\circ\eta} (A,C)_\psi$ is
Galois iff the extension $B\hookrightarrow A$ is Galois 
(with the
canonical entwining map $\psi$).}
\eex

\begin{theorem}
\label{thm.ind}
Let $(\alpha,\gamma)$ measure  $(\tA,\tC)_\tpsi$
to $ (A,C)_\psi$. Then the
following are equivalent:

(1) The functors $-\widehat{\square}_C\tC$, $-\hat{\tens}_\tA A$
are inverse equivalences.

(2) The functors $-\widehat{\square}_C\tC$, $-\hat{\tens}_\tA A$
are exact  and $(\tA,\tC)_\tpsi\mea\alpha\gamma (A,C)_\psi$ is Galois.

(3) The functor $-\hat{\tens}_\tA A$ is faithfully exact 
and $(\tA,\tC)_\tpsi\mea\alpha\gamma (A,C)_\psi$ is Galois.

(4) The functor $-\widehat{\square}_C\tC$
is faithfully exact and 
$(\tA,\tC)_\tpsi\mea\alpha\gamma (A,C)_\psi$ is Galois.
\end{theorem}
\proof
Recall that a functor is exact (resp. faithfully exact) if it preserves
(resp. preserves and reflects) exact sequences. 
(1) clearly implies (2), (3) and (4). To show
that (2) implies (1), first notice that $(A\otimes
C)\widehat{\square}_C\tC \cong A\otimes \tC$ in ${}_A\CM^\tC$ via
$A\otimes \eps \tens \tC$ and $(\mu\otimes C\otimes\tC)\circ(A\otimes \gamma \otimes
C)\circ(A\otimes \tilde{\Delta})$. Let $\tilde{M} = (A\otimes
C)\widehat{\square}_C\tC$. Then for all $M\in \M_A^C(\psi)$, 
$M\otimes_A\tilde{M}$ is in $\M_\tA^\tC(\tpsi)$ via
$M\otimes_A\mu_{\tilde{M}}$, $M\otimes_A\Delta_{\tilde{M}}$ and there is
a commutative diagram:
$$
\begin{CD}
(M\otimes_A\tilde{M})\otimes\tA\otimes A @>>> (M\otimes_A\tilde{M})\tens A 
@>{\hat{\pi}_{M\otimes\tilde{M}}}>>
(M\otimes_A\tilde{M})\hat{\tens}_\tA A @>>> 0\\
@VV{\cong}V                          @VV{\cong}V       @VVfV \\
M\otimes_A(\tilde{M}\otimes\tA\otimes A) @>>> M\otimes_A(\tilde{M}\tens A) 
@>{M\otimes_A\hat{\pi}_{\tilde{M}}}>>
M\otimes_A(\tilde{M}\hat{\tens}_\tA A) @>>> 0
\end{CD}
$$
The top row is exact since it is the defining sequence of $\hat\otimes$.
The bottom row is the defining sequence of $\hat\otimes$ tensored with
$M$ and thus is exact since the tensor product is right exact. 
Therefore the map $f$ (constructed from the diagram) is an 
isomorphism, and we have:
$$
(M\tens\tC)\hat{\tens}_\tA A\cong (M\otimes _A((A\otimes
C)\widehat{\square}_C \tC))\hat{\tens}_\tA A\cong M\otimes _A(((A\otimes
C)\widehat{\square}_C \tC)\hat{\tens}_\tA A).
$$
Thus we can consider the following  commutative diagram
\begin{equation}
\begin{CD}  
0 @>>>   (M\widehat{\square}_C\tC)\hat{\otimes}_{\tilde{A}}A  @>>>  
(M\otimes \tC)\hat{\otimes}_\tA A 
  @>>>( M\otimes C\otimes \tC)\hat{\otimes}_\tA{A} \\  
  @.      @VV\Phi_{M}V  
@VV{M}\otimes_{A}\Phi_{A\otimes C}V 
@VV({M}\otimes{C})\otimes_{{A}}\Phi_{{A}\otimes 
{C}}V\\  
0 @>>>  {M}   @>\Delta_{M}>>   {M} \otimes {C}  
  @>\ell_{{M}{C}}>> {M}\tens {C}\otimes {C}  
\end{CD}    
\label{1cd}
\end{equation}
The top row is the defining sequence of $M\widehat{\square}_C\tC$ acted
upon by $-\hat{\otimes}_\tA A$  and thus is exact by the exactness of   
$-\hat{\otimes}_\tA A$.
In the bottom row, $\ell_{{M}{C}} =
\Delta_{{M}}\otimes {C} - {M}\otimes \Delta$ and hence
the sequence  is exact by the definition of the coproduct. Since
$(\alpha,\gamma)$ is a Galois measuring, the maps in the second 
and the third
columns are bijective and thus so is $\Phi_{{M}}$. 

Now, reversing the arrows in the above diagram, interchanging 
$\widehat{\square}$ with $\hat{\tens}$, $A$ with $C$, 
tilded expressions with the
non-tilded ones, coactions with actions, and $\Phi$ with $\Psi$ 
one obtains the diagram from which
one deduces that also $\Psi_{\tilde{M}}$ is bijective provided the functor
$-\widehat{\square}_C\tC$ is exact. 

Next we show that (3) implies (2). From (\ref{1cd}) we know that
$\Phi_{{M}}$ is bijective for any ${M}\in
\M^{{C}}_{{A}}({\psi})$. Therefore for any exact
sequence ${M}_1\to{M}_2\to{M}_3$ of objects in
$\M^{{C}}_{{A}}({\psi})$ there is an exact sequence
$
({M}_1\widehat{\square}_{{C}}\tC)\hat{\otimes}_\tA{A} \to
({M}_2\widehat{\square}_{{C}}\tC)\hat{\otimes}_\tA{A}\to
({M}_3\widehat{\square}_{{C}}\tC)\hat{\otimes}_\tA{A}.
$
Since $-\hat{\tens}_\tA A$ reflects exact sequences,  
there is an exact sequence 
$
{M}_1\widehat{\square}_{{C}}\tC\to
{M}_2\widehat{\square}_{{C}}\tC\to
{M}_3\widehat{\square}_{{C}}\tC,
$
i.e.,  $-\widehat{\square}_{{C}}\tC$ is exact as required.  
Similarly one shows that (4) implies
(2). \endproof

Theorem~\ref{thm.ind} applied to entwining structures of
Example~\ref{ex.cat}(3) and  measurings coming from morphisms of
entwining structures gives \cite[Theorem~2.8]{CaeRai:ind} 
(and Proposition~\ref{prop.adjoint.n} gives 
\cite[Theorem~1.3]{CaeRai:ind}). Furthermore we obtain the
following  generalisation of  \cite[Theorem~3.7]{Sch:pri}
  
\bco\label{equivalence}  
For an entwining structure $(A,C)_\psi$ the following are  
equivalent:  
  
(1) $A(B)^C$ is a $C$-Galois extension with the canonical entwining  
map $\psi$ and $A$ is faithfully flat as a left $B$-module.  
  
(2) $A\in \M_A^C(\psi)$ and the functor $\M^C_A(\psi) \to \M_B$,  
$M\mapsto M_0$ is an equivalence.  
\eco  
\proof A left $B$-module $A$ is faithfully  
flat iff the functor $-\otimes_BA$ is faithfully exact, hence
the assertion follows  by
applying Theorem~\ref{thm.ind} to Example~\ref{just}. \endproof
 
\bco\label{cor.equ.cleft}   
If $(A,C)_\psi$ is the canonical entwining structure associated to a cleft 
 $C$-Galois extension $A(B)^C$ then 
$\M_A^C(\psi)$ is equivalent to $\M_B$.
\eco

\proof $A\cong   
B\tens C$ as objects in ${}_B\CM^C$, so $A$ is a faithfully flat
left $B$-module. \endproof

The following proposition, which is an $(A,C)_\psi$-module version of  
\cite[Theorem~2.11]{DoiTak:miy}, gives a criterion for  faithful  
flatness of a $C$-Galois extension $A(B)^C$  
\bpr\label{int}  
Let $A(B)^C$ be a $C$-Galois extension and assume that there exists a  
linear map $\phi:C\to A$ such that $1\sw 0\phi(1\sw 1) = 1$ and  
$$  
\psi(c\sw 1\otimes \phi(c\sw 2)) = \phi(c) 1\sw 0\otimes 1\sw 1, \qquad  
\forall c\in C.  
$$  
If either $A$ is flat as a left $B$-module or for all $b\in B$ and  
$c\in C$, $b_\alpha\phi(c^\alpha) = \phi(c)b$ then $A$ is faithfully  
flat as a left $B$-module.  
\epr  
\proof We are in the setting of  Example~\ref{just}, thus it suffices to 
show that the  
functor $M\mapsto M_0$ is  an equivalence and then use 
Corollary~\ref{equivalence} to deduce the assertion. First notice that for all $a\in A$, we have 
$a\sw 0\phi(a\sw 1)  
\in B$. For any right $B$ module $V$ the adjunction $\Psi_V: V\to
(V\otimes_BA)_0$ is simply $v\mapsto v\otimes_B 1$ and has the inverse
$\Psi^{-1}_V : (V\otimes_  
BA)_0 \to V$, $\sum_i v^i\otimes_Ba^i \mapsto \sum_iv^i\cdot (a^i\sw  
0\phi(a^i\sw 1))$. Now consider 
the commutative diagram (\ref{1cd}) for the
measuring  of Example~\ref{just}. If $A$ is  
a flat left  $B$-module then the top sequence is exact and thus  
$\Phi_M$ is bijective. Therefore the equivalence of categories holds  
in this case. On the other hand we have an exact sequence  
$$  
0\to \overline{M\otimes_B A} \to M\otimes_B A \to (M\otimes C)\otimes_B
A.   
$$  
The elements of $\overline{M\otimes_B A}$ are characterised by the property  
$\sum_i \Delta_M(m^i)\otimes_B a^i = \sum_i m^i\cdot 1\sw 0 \otimes 1\sw 1
\otimes_B  
a^i$. $M_0\otimes_B A$ is included in $\overline{M\otimes_B A}$ 
canonically. Consider the map $M\to M_0$, $m\mapsto m\sw  
0\cdot \phi(m\sw 1)$. If for all $b\in B$ and  
$c\in C$ $b_\alpha\phi(c^\alpha) = \phi(c)b$ then this map is a right  
$B$-module homomorphism. This implies that the map $\overline{M\otimes_B  
A}\to M_0\otimes_B A$, $\sum_i m^i\otimes_B a^i \mapsto \sum_i  
m^i\sw 0 \cdot \phi(m^i\sw 1) \otimes_B a^i$ is well-defined. It is an  
easy exercise to verify that this map is an inverse of the canonical   
inclusion $M_0\otimes_B A \hookrightarrow \overline{M\otimes_B A}$. Thus we  
conclude that also in this case the top sequence in (\ref{1cd}) (for a
measuring of Example~\ref{just}) is  
exact so that the functor $M\mapsto M_0$ is an equivalence. \endproof

\section*{\normalsize\sc\centering 4. Associated  
modules $A\square_CV$}\vspace{-.6\baselineskip}   
\setcounter{section}{4}     
\setcounter{definition}{0}  
  
In this section we construct the left 
$B$-module $E$ for each $C$-Galois extension $A$ of $B$ and a left 
$C$-comodule $V$, and then study the properties of $E$. This construction
is a very general algebraic dualisation of  associating of a 
fibre bundle to a principal bundle. 
  
First recall the definition of a cotensor product \cite{MilMoo:str}.     
Let $C$ be a coalgebra and $V_R\in \M^C$, $V_L\in {}^C\M$.     
The cotensor product $V_R\square_{C}V_L$ is defined by the exact sequence    
\[    
0\longrightarrow V_R\square_{C}V_L\hookrightarrow V_R\ot V_L\st{\ell_{V_RV_L}}    
{\longrightarrow}V_R\ot C\ot V_L,    
\]    
where $\ell_{V_RV_L}$ is the coaction equalising map     
$\ell_{V_RV_L}=\Delta_{V_R}\ot V_L-V_R\ot\, _{V_L}\Delta$.    
    
\bde\label{vc-ext}     
Let $A(B)^C$ be a $C$-Galois extension. A left $B$-module $E$ is called 
a {\em left module associated to $A(B)^C$} iff there exists a left     
$C$-comodule $V$ such that  $E=A\square_C V$. In this case $E$ is denoted 
by $E(A(B)^C;V)$.
\ede     
Since $A(B)^C$ can be viewed as an object dual to a (generalised) 
principal bundle and $V$ is dual to a representation of a ``structure
group", $E$ can be viewed as an object dual to a fibre bundle associated
to a principal bundle. 
In particular, \cite[Lemma~3.1(i)]{Sch:pri} 
implies that a {\em quantum fibre bundle} of
\cite[Definition~A.3]{BrzMa:gau} associated to a Hopf-Galois extension
$A(B)^H$ is a left module associated to $A(B)^H$ provided the 
antipode in $H$ is bijective. As
should be expected, $E(A(B)^C;C) = A(B)^C$, since  
$A\cong A\square_CC$ via $\Delta_A:A\to A\square_CC$ and 
$A\otimes\eps:A\square_CC\to
A$  (cf. \cite[Lemma~2.2*]{Gug:ext}). Furthermore, if $A(B)^C$ is
a cleft $C$-Galois extension then $A\cong B\tens C$   
in ${}_B\CM^C$, and for any left $C$-comodule $V$ we have
 $E = (B\otimes C)\square_C V \cong B\otimes V$ in
${}_B\M$. This last statement reflects the classical fact that every
fibre bundle associated to a trivial principal bundle is trivial.
  
\bde\label{def.cross-section}     
Let $E$ be a left module
associated  to $A(B)^C$. Any left $B$-module map $s: E\to B$
is called a {\em cross-section} of $E$.   
\ede     
The space $~{\rm Hom}_{B-}(E,B)$  of all cross-sections   
of $E(A(B)^C;V)$ has a natural right $B$-module structure   
given by $(s\cdot b)(x) = s(x)b$, for all $b\in B$, $x\in E$. Let    
for a given $C$-Galois extension $A(B)^C$ and a left $C$-comodule $V$,   
${\rm Hom}_{\psi}(V, A)$ denote the space of all linear maps $V\to A$ such   
that for all  
$v\in V$  
\begin{equation}\label{con.phi}     
\psi(v\sw{-1}\tens \phi(v\sw 0)) = \phi(v)\Delta_A(1),  
\end{equation}     
where $\psi:C\otimes A\to A\otimes C$ is the canonical entwining map    
associated to $A(B)^C$. For all $\phi\in {\rm Hom}_{\psi}(V, A)$,   
$b\in B$ and $v\in V$ we have   
\begin{eqnarray*}   
\psi(v\sw{-1}\tens \phi(v\sw 0) b) & = &  \phi(v\sw 0)_\alpha \psi(
v\sw{-1}^{\alpha}\tens b) \qquad \qquad \mbox{\rm (by (\ref{diag.A}))}\\    
& = &\phi(v)1\sw 0\psi(1\sw 1\otimes b) \qquad \qquad \;\;\;\mbox{\rm
($\phi\in {\rm Hom}_\psi (V, A)$)} \\  
& = & (\phi(v)b)1\sw 0\tens 1\sw 1.
\qquad \qquad \;\;\;\; \mbox{\rm ($b\in B$)} 
\end{eqnarray*}   
This implies that ${\rm Hom}_{\psi}(V, A)$ is a right $B$-module with the action given   
by $(\phi\cdot b)(v) = \phi(v)b$. The following theorem reflects the classical equivalence
between cross-sections of a fibre bundle and  
equivariant functions on it (cf. \cite[4.8.1]{Hus:fib}).   
\begin{theorem}\label{pro.cross-section}     
Let $E = E(A(B)^C;V)$. If either $A$ is 
a flat right $B$-module or else $V$ is a coflat left $C$-comodule, then 
the right $B$-modules ${\rm Hom}_{B-}(E,B)$ and ${\rm Hom}_\psi     
(V,A)$ are isomorphic to each other. 
\end{theorem}     
    
\proof  The flatness (coflatness) assumption implies that
$(A\otimes_BA)\square_CV \cong
A\otimes_B(A\square_CV)$, canonically (cf. \cite[p.~172]{Sch:pri}). 
Thus there is a left
$A$-module isomorphism $\rho: A\otimes_B E \to A\otimes V$, 
obtained as a composition of $can\square_C V$ with the canonical
isomorphism $A\otimes C\square_CV \stackrel{\sim}{\to} A\otimes V$, i.e.,
$\rho = \mu\otimes V$, $\rho^{-1} = (can^{-1}\otimes
V)\circ(A\otimes{}_V\Delta)$. 
Following \cite{DoiTak:miy}, apply  ${\rm Hom}_{A-}(-,A)$ to $\rho$ to
deduce the isomorphism
$$
\theta :{\rm Hom}(V,A) \stackrel{\sim}{\to}{\rm Hom}_{B-}(E,A), \qquad
\theta(\phi) (\sum_ia^i\tens v^i)=   
\sum_ia^i\phi(v^i). 
$$
Notice that $\theta$ is a right $B$-module map. 
For any $\phi\in {\rm Hom}(V,A)$, $x=\sum_ia^i\tens v^i \in E$
$$
\Delta_A(\theta(\phi)(x)) = \sum_i\Delta_A(a^i\phi(v^i))=    
\sum_i a^i\sw 0\psi(a^i\sw 1\tens\phi(v^i)) 
= \sum_i a^i\psi(v^i\sw{-1}\tens\phi(v^i\sw{0})).
$$
Therefore $\theta(\phi)\in {\rm Hom}_{B-}(E,B)$ iff
\begin{equation}
\sum_i a^i\psi(v^i\sw{-1}\tens\phi(v^i\sw{0})) = \sum_i  
a^i\phi(v^i)1\sw 0\tens 1\sw 1,
\label{sec.1}
\end{equation}
since $B=A_0$ by Example~\ref{just}. 
Clearly, (\ref{con.phi}) implies (\ref{sec.1}). Applying (\ref{sec.1})
to $\rho^{-1}(1\otimes v)$
one easily finds that
(\ref{sec.1}) implies (\ref{con.phi}). Therefore $\theta$ restricts to
${\rm Hom}_\psi(V,A) \stackrel{\sim}{\to} {\rm Hom}_{B-}(E,B)$ as a
right $B$-module map. \endproof

Viewing a $C$-Galois extension as a left module associated to itself, 
 one can 
state Proposition~\ref{int} as follows   
\bpr\label{crit.ff}  
If a  $C$-Galois extension $A(B)^C$ admits a unital $B$-bimodule map  
$s :A\to B$ then $A$ is faithfully flat as a left $B$-module.  
\epr  
\proof   
We view $A(B)^C$ as $E(A(B)^C;C)$. Then $\phi \in {\rm   
Hom}_{\psi}(C,A)$ iff for all $c\in C$, $\psi(c\sw 1\tens \phi(c\sw   
2)) = \phi(c)\Delta_A(1)$. The cross-sections are simply  
left $B$-module maps $A\to B$. Since $C$ is coflat in ${}^C\M$, by 
Theorem~\ref{pro.cross-section}, there is an 
isomorphism of  right  
$B$-modules $\theta: {\rm Hom}_\psi     
(C,A)\stackrel{\sim}{\to}  {\rm Hom}_{B-}(A,B)$, 
$\theta(\phi)(a) = a\sw 0\phi(a\sw 1)$ and  
$\theta^{-1}(s)(c) =c\su 1s(c\su 2)$, where $c\su 1\otimes_B c\su 2 =
can^{-1}(1\otimes c)$. We now assume that there is a unital  
$B$-bimodule map $s: A\to B$, and let $\phi = \theta^{-1}(s)$. 
Since $s$ is unital  $1\sw 0\phi(1\sw 1)  
:= 1\sw 0\theta^{-1}(s)(1\sw 1) = 1\sw 01\sw 1\su 1 s(1\sw 1\su 2) =  
s(1) =1$. Furthermore for any $b\in B$ and $c\in C$   
\begin{eqnarray*}  
b_\alpha \phi(c^\alpha) & = & c\su 1(c\su 2b)\sw 0 \phi ((c\su 2b)\sw  
1) = c\su 1(c\su 2b)\sw 0 (c\su 2b)\sw 1\su 1s((c\su 2b)\sw 1\su 2)\\  
& = & c\su 1 s(c\su 2 b) =\phi(c) b,  
\end{eqnarray*}  
where we used the definition of the canonical entwining structure to  
derive the first equality, then the following property of the
translation map  (cf. \cite[Remark~3.4]{Sch:rep}) 
\begin{equation}
a\sw 0 a\sw 1\su 1\otimes_B a\sw 1\su 2 = 1\otimes_B a, \qquad \forall
a\in A,
\label{tra.2}
\end{equation}
 to derive the third one, 
and the right $B$-module property of $s$ to obtain the last  
equality. Thus we conclude that $\phi$ satisfies all the assumptions  
of Proposition~\ref{int} and hence the assertion follows. \endproof

Next we establish  the equivalence between a certain class of  
cross-sections of $A(B)^C$ and cleaving maps (cf. \cite{MasDoi:gen} for
an interesting special case).  This reflects the classical equivalence
between cross-sections and trivialisations of a principal bundle (cf.
\cite[4.8.3]{Hus:fib}) 
\begin{prop}\label{pro.sec.cleft}   
A $C$-Galois extension $A(B)^C$ is cleft if and only if there exists   
a  cross-section $s\in{}{\rm Hom}_{B-}(A,B)$   
 such that $\hat{s} :=   
(s\tens C)\circ\Delta_A :A\to B\tens C$ is a bijection.  
\end{prop}   
\proof  
Assume first that $A(B)^C$ is cleft with a cleaving map $\Phi:C\to   
A$. Then (\ref{phi-1}) implies
that $\Phi^{-1}\in {\rm   
Hom}_{\psi}(C,A)$. By Theorem~\ref{pro.cross-section} there is a   
cross-section $s=\theta(\Phi^{-1})$. Explicitly $s(a) = a\sw   
0\Phi^{-1}(a\sw 1)$. The induced map $\hat{s}$ is thus $\hat{s}(a) = a\sw   
0\Phi^{-1}(a\sw 1)\tens a\sw 2$ and has the inverse $b\tens c \mapsto   
b\Phi(c)$ as in Proposition~\ref{ex.cleft}.   
   
Assume now that there exists $s\in {\rm Hom}_{B-}(A,B)$ such that the   
map $\hat{s}$ is bijective. Since $s$ is a left $B$-module map,   
$\hat{s}$ is a morphism in ${}_B\CM^C$, where $B\tens C$ is viewed   
as an object in ${}_B\CM^C$ via $\mu\tens C$ and $B\tens\Delta$. This  
implies that also    
$\hat{s}^{-1}$ is a morphism in ${}_B\CM^C$. Note also that $s =   
(B\tens \eps)\circ \hat{s}$. Using Theorem~\ref{pro.cross-section}   
we consider $\tilde{\Phi}\in {\rm   
Hom}_{\psi}(C,A)$ given by $\tilde{\Phi} = \theta^{-1}(s)$, and also a   
map $\Phi :C\to A$, $\Phi: c\mapsto \hat{s}^{-1}(1\tens c)$. We will   
show that $\Phi$ and $\tilde{\Phi}$ are convolution inverses to each   
other. For any $c\in C$ one has   
\begin{eqnarray*}  
\Phi(c\sw 1)\tilde{\Phi}(c\sw 2) & = &\hat{s}^{-1}(1\tens c\sw   
1)\theta^{-1}(s)(c\sw 2) = \hat{s}^{-1}(1\tens c\sw   
1)c\sw 2\su 1s(c\sw 2\su 2)\\
&=& \hat{s}^{-1}(1\tens c)\sw 0   
\hat{s}^{-1}(1\tens c)\sw 1 \su 1s(\hat{s}^{-1}(1\tens c)\sw 1\su 2)
\qquad~\;\mbox{\rm ($\hat{s}$ is $C$-colinear)}\\
& = & s(\hat{s}^{-1}(1\tens c)) \qquad~\qquad~\qquad~\qquad~\qquad~\qquad~\qquad \mbox{\rm
(by (\ref{tra.2}))}\\
& = &  
((B\tens\eps)\circ\hat{s}\circ\hat{s}^{-1})(1\tens c) = \eps(c),  
\end{eqnarray*}
where $c\su 1\otimes_Bc\su 2 = can^{-1}(1\otimes c)$ as before. 
On the other hand, for all $a\in A$, one finds   
$$
a\sw 0\tilde{\Phi}(a\sw 1)\Phi(a\sw 2) =  a\sw 0 \theta^{-1}(s)(a\sw   
1)\hat{s}^{-1}(1\tens a\sw 2)  
 = a\sw 0 a\sw 1\su 1s(a\sw 1\su 2)\hat{s}^{-1}(1\tens a\sw 2).   
$$
Making use of (\ref{tra.2}) and   the fact that   
$\hat{s}^{-1}$ is a left $B$-module map one concludes   
$$   
a\sw 0\tilde{\Phi}(a\sw 1)\Phi(a\sw 2) =  s(a\sw 0) \hat{s}^{-1}(1\tens a\sw   
1) = \hat{s}^{-1}(s(a\sw 0)\tens a\sw 1) = \hat{s}^{-1}(\hat{s}(a)) =   
a.   
$$   
Finally, using the definition of the translation map one obtains for all   
$c\in C$   
$$   
\tilde{\Phi}(c\sw 1)\Phi(c\sw 2) = c\su 1 c\su 2\sw 0\tilde{\Phi}(c\su   
2\sw 1)\Phi(c\su 2\sw 2) = c\su 1 c\su 2 = \eps(c).   
$$    
{}From the fact that $\tilde{\Phi}\in {\rm   
Hom}_{\psi}(C,A)$ it is clear that its convolution inverse $\Phi$ is   
a right $C$-comodule map. Thus we have proven that $A(B)^C$ is   
cleft. \endproof

Since $\hat{s}$ is a  morphism in ${}_B\CM^C$, 
Proposition~\ref{pro.sec.cleft} states that a 
$C$-Galois extension $A(B)^C$ is cleft if and only if $A\cong   
B\tens C$ as objects in ${}_B\CM^C$ (cf. \cite[Theorem~9]{DoiTak:cle})   

\section*{\normalsize\sc\centering 5. Associated  
modules $(V\otimes A)_0$}\vspace{-.6\baselineskip}   
\setcounter{section}{5}     
\setcounter{definition}{0}  
  
In this section we construct the right 
$B$-module $\bar{E}$ for each $C$-Galois extension $A$ of $B$ and a right 
$C$-comodule $V$, and then study  properties of $\bar{E}$. This construction
is another  algebraic dualisation of  associating of a 
fibre bundle to a principal bundle.

\bde\label{vr}   
Let $A(B)^C$ be a coalgebra Galois extension. A  right $B$-module 
$\bar{E}$ is called  
a {\em right module associated to $A(B)^C$} iff there exists a   
right $C$-comodule $V$ such that $\bar{E} = (V\otimes A)_0$ 
(cf. Example~\ref{just}), where   
$V\otimes A$ is an $(A,C)_\psi$-module of the canonical
entwining structure as in  
Corollary~\ref{MC1}(2). In this case $\bar{E}$ is denoted by
$\bar{E}(V;A(B)^C)$.
\ede  
  
The right $B$-module $\bar{E}$ consists of all elements $\sum_iv^i\otimes  
a^i \in V\otimes A$ such that $\sum_i v^i\sw 0\otimes \psi(v^i\sw 1 \otimes  
a^i) = \sum_i v^i\otimes a^i1\sw 0\otimes 1\sw 1$. In the case of a   
Hopf-Galois extension, if $V$ is a right comodule algebra the above definition  
of $\bar{E}$ coincides with the definition of a quantum associated fibre  
bundle in \cite[Definition~A.1]{Haj:str}.  Thus, similarly as in Section~5,
the constructed right module $\bar{E}$ can be viewed geometrically as an
object dual to a fibre bundle associated to a principal bundle of which
$A(B)^C$ is a dual object.
  
\bpr\label{triv.r}  
If $A(B)^C$ is cleft then for all right comodules $V$, $\bar{E} =
(V\otimes A)_0$ is isomorphic to $V\otimes B$ as a right $B$-module.  
\epr  
\proof The isomorphism and its inverse are:  
$$  
V\otimes B \to \bar{E}, \qquad v\otimes b \mapsto v\sw 0 \otimes \Phi^{-1}  
(v\sw 1) b,  
$$  
$$  
\bar{E}\to V\otimes B, \qquad  \sum_i v^i\otimes a^i \mapsto v^i\sw 0 \otimes  
\Phi(v^i\sw 1)a^i,  
$$  
where $\Phi$ is a cleaving map. To see that the output of the first of these maps is in $\bar{E}$ we  
compute  
\begin{eqnarray*}  
v\sw 0\otimes \psi(v\sw 1\otimes \Phi^{-1}(v\sw 2) b) & =  &
v\sw 0\otimes\Phi^{-1}(v\sw 1)1\sw 0\psi(1\sw 1\otimes b) \qquad~\;\;\;
\mbox{\rm (by (\ref{diag.A}) and (\ref{phi-1}))}\\  
& = & v\sw 0\otimes\Phi^{-1}(v\sw 1)b1\sw 0\otimes 1\sw 1. 
\end{eqnarray*}  
To verify
that the output   
of the second of the above maps is in $V\otimes B$ we use the fact  
that $A$ is an $(A,C)_\psi$-module and that $\Phi$ is a right
$C$-comodule map to compute  
$$
\sum_i v^i\sw 0 \Delta_A(\Phi(v^i\sw 1)a^i) =  \sum_i v^i\sw 0
\Phi(v^i\sw 1)   
\psi(v^i\sw 2\otimes a^i)  
= \sum_i v^i\sw 0 \Phi(v^i\sw 1)a^i1\sw 0\otimes  
1\sw 1. 
$$ 
It is obvious that the above maps are inverses to each other and that they  
are right $B$-module homomorphisms. \endproof

\bde\label{sec.r}  
Let $\bar{E}$ be a right module associated to $A(B)^C$. 
Any right $B$-module map $s: \bar{E}\to B$ is called a
{\em cross-section}  
of $\bar{E}$.
\ede  
  
The space ${\rm Hom}_{-B}(\bar{E},B)$ of cross-sections of 
$\bar{E}$ has a   
natural left $B$-module structure.  
Since $\Delta_A$ is left linear over $B$, the space ${\rm 
Hom}^{-C}(V,A)$ of all right $C$-comodule maps $\phi:V\to A$ has a 
left $B$-module structure $b\cdot \phi:    
v\mapsto b\phi(v)$.  
  
\bth\label{cross.r}  
For any $\bar{E}(V;A(B)^C)$,  
if $A$ is faithfully flat as a left $B$-module then   
${\rm Hom}^{-C}(V,A)\cong {\rm Hom}_{-B}(\bar{E},B)$ as left   
$B$-modules.  
\ethe  
\proof  
By Corollary~\ref{adjoint.cor}(2) there is a natural isomorphism $\zeta_{A,V}:{\rm 
Hom}_{-A}^{-C}(V\otimes A,A) \stackrel\sim\to {\rm Hom}^{-C}(V,A)$, 
given by $\zeta_{A,V}(s)(v)=s(v  
\otimes 1)$, $\zeta^{-1}_{A,V}(\phi)(v\otimes a) = \phi(v)a$.
It is an easy exercise to verify 
that $\zeta_{A,V}$ preserves the left $B$-module structure. 
If $A$ is faithfully flat as a left $B$-module  
then by Corollary~\ref{equivalence} the functor $\M_A^C(\psi)\to \M_B$,  
$M\mapsto M_0$ is an equivalence. Therefore  the space 
${\rm Hom}_{-A}^{-C}(V\otimes A,A)$ is isomorphic to the space 
of morphisms $(V\otimes A)_0\to B$ in $\M_B$. Since the latter is 
precisely  
${\rm Hom}_{-B}(\bar{E},B)$ and the equivalence preserves the left 
$B$-module structure introduced on the spaces of morphisms, we conclude 
that there is an isomorphism of left $B$-modules  
${\rm Hom}_{-B}(\bar{E},B)\cong {\rm Hom}^{-C}(V,A)$.
\endproof  
  
\section*{\normalsize\sc\centering 6. Bijectivity of $\psi$ and the
relationship  
between $E$ and $\bar{E}$}\vspace{-.6\baselineskip}  
\setcounter{section}{6}  
\setcounter{definition}{0}  
 
In this section we study a relationship between left and right 
modules associated to a $C$-Galois extension. In particular we show that
 if $V$ is finite 
dimensional then $E(A(B)^C;V)$ can be 
identified with the module of cross-sections of $\bar{E}(V^*;A(B)^C)$, and vice 
versa. Also we show that if the canonical entwining 
map $\psi$ is bijective then $A\square_C V$ can be identified with the 
space of coinvariants $_0(A\otimes V)$ of the left coaction of $C$ on 
$A\otimes V$. The module $_0(A\otimes V)$ can be then viewed as 
$\bar{E}$ for the left $C$-Galois extension $A$ of $B$. 
 
If $V$ is a finite-dimensional left $C$-comodule, the dual vector space 
$V^*$ is viewed as a right $C$-comodule via 
$
\langle v^*\sw 0, v\rangle v^*\sw 1 = v\sw{-1}\langle v^*,v\sw 0 
\rangle, 
$
for all $v\in V$, $v^*\in V^*$, where $\langle\cdot,\cdot\rangle :V^*\otimes V\to k$ denotes the 
non-degenerate pairing.
 
\bpr\label{duality} 
Let $A(B)^C$ be a $C$-Galois extension, $V$ be a finite 
dimensional left 
$C$-comodule, and $V^*$ the dual right $C$-comodule. Then: 
 
(1) The left $B$-modules $E=A\square_C V$ and ${\rm Hom}^{-C}(V^*, A)$ 
are isomorphic to each other. 
 
(2) The right $B$-modules $\bar{E} = (V^*\otimes A)_0$ and ${\rm 
Hom}_{\psi}(V, A)$ are isomorphic to each other. 
\epr 
\proof (1) It is well-known that the vector spaces 
${\rm Hom}^{-C}(V^*,A)$ and $A\square_C V$ are isomorphic to each other with 
the isomorphism $\phi\mapsto \sum_i a^i\otimes v^i$, such that for all 
$v^*\in V^*$, $\phi(v^*) = \sum_ia^i\langle v^*,v^i\rangle$. Clearly,
this is also an isomorphism of left $B$-modules.

(2) We identify $\phi\in {\rm Hom}(V, A)$ with $\sum_iv^{*i}\otimes 
a^i \in V^*\otimes A$ via $\phi(v) = \sum_i\langle v^{*i},v\rangle 
a^i$. Clearly this identification is an isomorphism of right 
$B$-modules. We have 
$$ 
\psi(v\sw{-1}\otimes \phi(v\sw 0)) = \sum_i\langle v^{*i},v\sw 0\rangle 
\psi(v\sw{-1}\otimes a^i) = \sum_i\langle v^{*i}\sw 0,v\rangle 
\psi(v^{*i}\sw{1}\otimes a^i). 
$$ 
On the other hand  
$ 
\phi(v)1\sw 0\otimes 1\sw 1 = \sum_i\langle v^{*i},v\rangle 
a^i1\sw 0\otimes 1\sw 1$ which implies that $\phi \in {\rm 
Hom}_\psi(V, A)$ if and only if $\sum_iv^{*i}\otimes 
a^i \in (V^*\otimes A)_0 = \bar{E}$. \endproof 
 
\bco\label{equiv.2} 
(1) Let $E = E(A(B)^C;V)$. If $V$ is finite dimensional and 
$A$ is faithfully flat as a 
left $B$-module then the left $B$-modules $E$ and 
${\rm Hom}_{-B}({\rm Hom}_\psi(V,A), B)$ are isomorphic to each other. 
 
(2)  Let $\bar{E} = \bar{E}(V;A(B)^C)$ with finite dimensional 
$V$. If either $A$ is flat as a 
right $B$-module or else  $V$ is coflat as a left $C$-comodule, then the right $B$-modules $\bar{E}$
and 
${\rm Hom}_{B-}({\rm Hom}^{-C}(V,A), B)$ are isomorphic to each other. 
\eco 
\proof (1) By Proposition~\ref{duality} one identifies $E$ with ${\rm 
Hom}^{-C}(V^*,A)$ and also ${\rm Hom}_\psi(V,A)$ with
$\bar{E}(V^*;A(B)^C)$, and then applies Theorem~\ref{cross.r} to deduce 
the isomorphism of left $B$-modules. 
 
(2) By Proposition~\ref{duality} one identifies $\bar{E}$ with ${\rm 
Hom}_{\psi}(V^*,A)$ which by Theorem~\ref{pro.cross-section} is
isomorphic to the right module of cross-sections ${\rm
Hom}_{B-}(A\square_CV^*,B)$. Then one applies  Proposition~\ref{duality}
again to deduce the required isomorphism.
\endproof 
 
\bre{\rm For any left module $E(A(B)^C;V)$ there is a left $B$-module map $E\to 
{\rm Hom}_{-B}({\rm Hom}_\psi(V,A), B)$ given by $\sum_i a^i\otimes 
v^i\mapsto s$ where $s(\phi) = \sum_ia^i\phi(v^i)$. Similarly for any 
$\bar{E}(V;A(B)^C)$ there is a right $B$-module map 
$\bar{E}\to {\rm Hom}_{B-}({\rm Hom}^{-C}(V,A), B)$ given by 
$\sum_iv^i\otimes a^i\mapsto s$, where $s(\phi) = 
\sum_i\phi(v^i)a^i$.} 
\ere 

The remaining part of this section is devoted to studies of the
relationship  between $E$ and $\bar{E}$ in the case when the canonical
entwining map is bijective. 

\bex\label{ex.psi-1}
{\rm  (1) Let $H$, $A$, $C$ and $\psi$ be as in Example~\ref{ex.cat}(3). If the
antipode $S$ in $H$ is bijective then $\psi$ is bijective. Explicitly
$\psi^{-1}(a\otimes c) = c\cdot S^{-1} a\sw 1\otimes a\sw 0$, $\forall
a\in A, c\in C$. 

(2) For a Hopf-Galois extension $A(B)^H$, the canonical entwining map
$\psi: h\otimes a\mapsto a\sw 0\otimes ha\sw 1$  
is bijective if and only if the antipode in $H$ is bijective.} 
\eex
\proof (1) is proven by a straightforward  computation. To prove (2) 
consider the linear map  $\psi_H: H\otimes H\to H\otimes H$,
$h\otimes h' \mapsto h'\sw 1\otimes hh'\sw 2$. It is well-known that
$\psi_H$ is bijective if and only if the antipode is bijective. Notice
that  $A\otimes_B\psi = (can^{-1}\otimes H)\circ (A\otimes \psi_H)\circ
(can\otimes H)$. This completes the proof. \endproof

\ble\label{lem.psi-1} 
Let $A(B)^C$ be a $C$-Galois extension and assume that the canonical 
entwining map $\psi$ is bijective. Then: 
 
(1) $A$ is a left $C$-comodule with the coaction $_A\Delta(a) = 
\psi^{-1}(a1\sw 0\otimes 1\sw 1)$. 
 
(2) The canonical right $A$-module left $C$-comodule map 
$can_L:A\otimes_B A\to C\otimes A$, $a\otimes a'\mapsto {}_A\Delta(a)a'$ 
is bijective. 
 
(3) The algebra $B$ is isomorphic to 
$$
\bar{B} :=  \{ b\in A \; | \; \forall a\in A,  
\; {}_A\Delta(a b) = _A\Delta(a)b\}
= \{b\in A \; |\; _A\Delta(b) = 
\psi^{-1}(1\sw 0\otimes 1\sw 1)b\}. 
$$
\ele 
\proof (1) is proven in the way analogous to the proof that 
$\psi$ induces a right $C$-coaction on $A$ via $\Delta_A(a) = 1\sw 
0\psi (1\sw 1\otimes a)$. To prove (2) notice that $can_L = 
\psi^{-1}\circ can$ thus it is well defined and a bijection. Next take 
$a\in A$ and $b\in B$, then $_A\Delta(ab) = \can_L(ab\otimes_B 1) = 
can_L(a\otimes_B b) = \can_L(a\otimes_B 1)b = {}_A\Delta(a) 
b$, i.e. ${b}\in \bar{B}$. On the other hand take 
$\bar{b}\in\bar{B}$. Then $_A\Delta(\bar{b}) = \psi^{-1}(1\sw 0 \otimes 1\sw 
1)\bar{b}$, i.e. $\psi^{-1}(\bar{b}1\sw 0\otimes 1\sw 1) = 
\psi^{-1}(1\sw 0 \otimes 1\sw 1)\bar{b}$. Applying $\psi$ one obtains 
$\bar{b}1\sw 0\otimes 1\sw 1 = 1\sw 0\psi(1\sw 1\otimes \bar{b}) = 
\Delta_A(\bar{b})$, i.e. $\bar{b}\in B$ by Example~\ref{just}. \endproof  
 
Having $\psi^{-1}:A\otimes C\to C\otimes A$ one can consider category 
$_A^C\M(\psi^{-1})$ the objects of which are left $A$-modules and left 
$C$-comodules $M$ such that 
$$ 
_M\Delta(a\cdot m) = \psi^{-1}(a\otimes m\sw{-1})\cdot m\sw 0. 
$$ 
The morphisms are left $A$-module left $C$-comodule maps. In 
particular, if the canonical
entwining map $\psi$ of $A(B)^C$ is bijective, then $A\in 
{}_A^C\M(\psi^{-1})$. One also defines a functor $_A^C\M(\psi^{-1})\to 
{}_B\M$, $M\mapsto {}_0M$, where $_0M := \{m\in M\; |\; _M\Delta(m) = 
\psi^{-1}(1\sw 0\otimes 1\sw 1)\cdot m\}$. If $V$ is a left 
$C$-comodule then $A\otimes V$ is an object in  $_A^C\M(\psi^{-1})$, 
where the left coaction is given by $_{A\otimes V}\Delta(a\otimes v) = 
\psi^{-1}(a\otimes v\sw{-1})\otimes v\sw 0$ and the action is
$\mu\otimes V$. Indeed, 
\begin{eqnarray*} 
_{A\otimes V}\Delta(a'a\otimes v) &= & \psi^{-1}(a'a\otimes 
v\sw{-1})\otimes v\sw 0  = v\sw{-1}_{\alpha\beta}\otimes 
a'^\beta a^\alpha\otimes v\sw 0 \\ 
& = & (a\otimes v)\sw{-1}_\beta \otimes a'^\beta\cdot (a\otimes v)\sw 0\\ 
& = & 
\psi^{-1}(a'\otimes (a\otimes v)\sw{-1})\cdot (a\otimes v)\sw 0, 
\end{eqnarray*} 
where  $\psi^{-1}(a\otimes c) = c_\alpha \otimes 
a^\alpha$. This implies that $_0(A\otimes V)$ is a left 
$B$-module. 
 
\bpr\label{relation} 
If the canonical
entwining map $\psi$ of $A(B)^C$ is bijective then:
 
(1) For any left $C$-comodule $V$, the left $B$-modules $_0(A\otimes 
V)$ and $A\square_C V$ are isomorphic to each other. 
 
(2) For any right $C$-comodule $V$, the right $B$-modules $(V\otimes 
A)_0$ and $V\square_C A$ are isomorphic to each other. 
\epr 
\proof (1) Take $\sum_ia^i\otimes v^i\in {}_0(A\otimes V)$. It means 
that 
$\sum_i\psi^{-1}(a^i\otimes v^i\sw{-1})\otimes v^i\sw 0 = \sum_i 
\psi^{-1}(1\sw 0\otimes 1\sw 1)a^i\otimes v^i$. Applying $\psi$ one 
obtains 
$$ 
\sum_ia^i\otimes v^i\sw{-1} \otimes v^i\sw 0 = \sum_i1\sw 0\psi(1\sw 
1\otimes a^i)\otimes v^i = \sum_ia^i\sw 0\otimes a^i\sw 1\otimes v^i. 
$$ 
Therefore $\sum_ia^i\otimes v^i \in A\square_CV$. To prove the second 
inclusion one repeats above steps in a reversed order. 
 
The proof of (2) is analogous to (1). \endproof 
 
Lemma~\ref{lem.psi-1} shows that any right $C$-Galois extension 
(i.e. with a right coaction) that has the bijective canonical 
entwining map, can 
be viewed equivalently as a left $C$-Galois extension (i.e. with a 
left coaction). Then Proposition~\ref{relation} yields that the module 
$E$ associated to a right $C$-Galois extension $A(B)^C$ plays the role of 
$\bar{E}$ when $A(B)^C$ is viewed as a left $C$-Galois 
extension. Similarly, the module $\bar{E}$ associated to the right 
extension corresponds to $E$ when $A(B)^C$ is viewed as a left $C$-Galois 
extension.

\section*{\normalsize\sc\centering Appendix. Dual  
results}\vspace{-.6\baselineskip}  
  
In this appendix we give dual version of the results described in  
Sections 2-6. Dual counterparts of statements 
given above are numbered with the same numbers decorated with  
stars. Proofs can be obtained by 
dualisation and thus are omitted.
\vspace{.5\baselineskip}  
  
\noindent { \normalsize\sc Definition~\ref{cge}* (\cite{BrzHaj:coa})} 
{\em Let $A$ be an algebra,  
$C$ a coalgebra and a right $A$-module with the action $\mu_C$,  and  
$B=C/I_C\,$,   
where $I_C\subset C$ is given by   
$$   
I_C = \span\{(c\cdot a)\sw 1a^*((c\cdot a)\sw 2)-c\sw 1a^*(c\sw 2\cdot a)
\;\;|\;\;  c\in C, a\in A, a^*\in A^*\}.  
$$  
We say  
that $C$ is an {\em algebra Galois coextension} (or $A$-Galois coextension)  
of $B$ iff the canonical left $C$-comodule right $A$-module   
map\vspace*{-2.5mm}  
\[  
cocan := (C\ot\mu_C)\circ(\Delta\tens A)\; :\;\; C\tens A \lra   
C\square\sb BC\vspace*{-2.5mm}  
\]  
is bijective. Here the coaction equalising map  
$\ell_{CC}$ is $\ell_{CC} = (C\otimes \pi_C)\circ\Delta -
(\pi_C\otimes C)\circ   
\Delta$, where $\pi_C :C\to B$ is the canonical surjection.  
Such an $A$-Galois coextension is denoted by  
$C(B)_A$.}\vspace{.5\baselineskip}   
  
We refer the reader to \cite[Section~3]{BrzHaj:coa}, where it is shown
that $B$ is a coalgebra, $\mu_C$ is left $B$-colinear and $cocan$ is
well-defined. Also in \cite{BrzHaj:coa} it is shown
 that every $A$-Galois      
coextension $C(B)_A$ induces a unique entwining map $\psi: C\tens A\to A\tens  
C$, $\psi =  
(\check{\tau}\tens C)\circ(C\tens\Delta)\circ cocan$,  
such that $C\in \M_A^C(\psi)$.    
Here $\check{\tau} : C\square_B C \to  
A$, $\check{\tau} := (\eps\tens C)\circ cocan^{-1}$ is the {\em  
cotranslation map}. This $\psi$ is called the {\em     
canonical entwining structure} associated to $C(B)_A$.   
  
A coextension $C(B)_A$ is {\em cocleft} if there exists a  
convolution invertible, right $A$-module map $\Phi:C\to A$
(cf. \cite[Definition~2.2]{MasDoi:gen}).  
The fact that $\Phi(c\cdot a) = \Phi(c)a$  
implies 
$$
\mu\circ(A\tens\Phi^{-1})\circ\psi = (\Phi^{-1}\tens\eps\circ\mu_C)\circ  
(\Delta\otimes A),  
$$
which, in turn, allows one to  
prove that $C\cong B\tens A$ as objects in ${}^B\CM_A$.  
  
For $C(B)_A$, ${\rm Aut}(C(B)_A)$ denotes the group of 
  left $B$-comodule, right $A$-module automorphisms of $C$ with the
product given by the composition of maps.\vspace{.5\baselineskip}  
  
\noindent{\normalsize\sc  Theorem~\ref{C-aut}*} {\em ${\rm
Aut}(C(B)_A)$  is isomorphic to the group $\CA(C)$ of 
convolution invertible maps $f:C\to A$ such that
\begin{spequation}{\ref{ad-cov}*}
\mu\circ(A\tens f)\circ\psi = \mu\circ (f\tens A),  
\end{spequation}
where $\psi$ is the canonical entwining map. The product in $\CA(C)$
is the convolution product.}\vspace{.5\baselineskip}  
  
\note{\proof Given $f:C\to A$ satisfying (\ref{ad-cov}*) define an  
$A$-Galois automorphism $F: C\to C$,  
$F :c\mapsto c\sw 1\cdot f(c\sw 2)$. Conversely, given an $A$-Galois  
automorphism $F:C\to C$ define $f:C\to A$, $f: c\mapsto  
\check{\tau}(c\sw 1, F(c\sw 2))$. \endproof }  
  
Notice that the condition (\ref{ad-cov}*) defining $\CA(C)$ can be
also understood as a twisted   
commutativity condition, since it explicitly reads for all $a\in A$, 
$c\in C$, $a_\alpha f(c^\alpha) = f(c)a$. 
If $C(B)_A$ is a cocleft $A$-Galois coextension, then $\CA(C)$
is  isomorphic to the group of convolution   
invertible maps $\gamma :B\to A$, since ${\rm End}_{-A}^{B-}(B\otimes A)
\cong {\rm Hom}(B,A)$ as algebras.\vspace{.5\baselineskip}   
  
\noindent {\normalsize\sc Example~\ref{just}*} Assume that  $C$ is an object in
$\M_A^C(\psi)$, and let $B$, $\pi_C$ be as in Definition~\ref{cge}*.
Then $(A,C)_\psi$ is measured to the trivial entwining structure 
$(k,B)_\sigma$ by $(\eps\circ\mu_C, \pi_C)$. With this measuring, for
all $V\in \M_k^B(\sigma) = \M^B$, $V\widehat{\square}_B^C = V\square_B
C$, while for all $M\in \M_A^C(\psi)$, $M\hat{\tens}_A k = M^0 :=
M/I_M$, where $I_M := \span\{m\cdot a - m\sw 0 \eps(m\sw 1\cdot a)\,|\; 
a\in A, \; m\in M\}$. Notice that $I_M =\span\{(m\cdot a)\sw 0 a^*((m\cdot a)\sw 1) - m\sw 0 
a^*(m\sw 1\cdot a)\,|\;  a\in A, \; m\in M, \; a^*\in A^*\}$. The
measuring $(\eps\circ\mu_C, \pi_C)$ is Galois iff the coextension
$C\twoheadrightarrow B$ is Galois. \vspace{.5\baselineskip}
 
\noindent{\normalsize\sc Corollary~\ref{equivalence}*} {\em   
For an entwining structure $(A,C)_\psi$ the following  are  
equivalent:  
  
(1) $C(B)_A$ is an $A$-Galois coextension with the canonical entwining  
map $\psi$ and $C$ is faithfully coflat as a left $B$-comodule (i.e. the
functor $-\square_BC$ is faithfully exact).  
  
(2) $C\in \M_A^C(\psi)$ and the functor $\M^C_A(\psi) \to \M^B$,  
$M\mapsto M^0$ is an equivalence.}\vspace{.5\baselineskip}  
 
\noindent{\normalsize\sc Corollary~\ref{cor.equ.cleft}*}{\em    
If $(A,C)_\psi$ is the canonical entwining structure associated to a cocleft 
 $A$-Galois coextension $C(B)_A$ then  
$\M_A^C(\psi)$ is equivalent to $\M^B$.}\vspace{.5\baselineskip}  
 
\noindent{\normalsize\sc Proposition~\ref{int}*} {\em Let $C(B)_A$ be an  
$A$-Galois coextension and assume that there exists a linear map  
$\phi:C\to A$ such that $\eps(c\sw 1\cdot\phi(c\sw 2)) = \eps(c)$ and  
$ 
a_\alpha\cdot\phi(c^\alpha) = \phi(c\sw 1)\eps(c\sw 2\cdot a).  
$ 
If either $C$ is coflat as a left $B$-comodule or for all $c\in C$, 
$\phi(c\sw 2)_\alpha\otimes \pi_C(c\sw 1^\alpha) = \phi(c\sw 1)\otimes 
\pi_C(c\sw 2)$, then $C$ is faithfully coflat as a left 
$B$-comodule.}\vspace{.5\baselineskip}  
 
\noindent {\normalsize\sc Definition~\ref{vc-ext}*} {\em Let $C(B)_A$ be  
an $A$-Galois coextension. A left $B$-comodule $E$ is called 
a {\em left comodule associated to $C(B)_A$} iff there exists a left     
$A$-module $V$ such that  $E=C\otimes_A V$. In this case $E$ is denoted 
by $E(C(B)_A;V)$.}\vspace{.5\baselineskip}    
 
For any $A$-Galois coextension, $E(C(B)_A;A) = C(B)_A$.
Also, if $C(B)_A$ is cocleft, then any $E(C(B)_A;V)$ 
is isomorphic to $B\tens V$ as a left   
$B$-comodule.\vspace{.5\baselineskip}   
  
\noindent{\normalsize\sc Definition~\ref{def.cross-section}*} {\em
Let $E$ be a left comodule associated      
to $C(B)_A$. Any left  
$B$-comodule map $s: B\to E$  is called a {\em cross-section} of 
$E$.}\vspace{.5\baselineskip}       
  
The space ${\rm Hom}^{B-}(B,E)$  of all cross-sections   
of $E(C(B)_A;V)$ has a natural right $B$-comodule
structure    
given by $\Delta_{{\rm Hom}^{B-}(B,E)}(s) = (s\tens B)\circ\Delta$. Let    
for a given $A$-Galois coextension $C(B)_A$ and a left $A$-module $V$,   
${\rm Hom}^{\psi}(C, V)$ denote the space of all linear maps $\phi:  
C\to V$ such  that   
\begin{spequation}{\ref {con.phi}*}  
{}_V\mu\circ(A\tens\phi)\circ\psi = 
(\phi\tens\eps\circ\mu_C)\circ(\Delta \otimes A),  
\end{spequation}  
where $\psi: C\tens A\to A\tens C$ is the canonical entwining map    
associated to      
$C(B)_A$. The space ${\rm Hom}^\psi
(C,V)$ is a right $C$-comodule via 
$(\phi\otimes \pi_C)\circ\Delta$.\vspace{.5\baselineskip}     
  
\noindent{\normalsize\sc Theorem~\ref{pro.cross-section}*} {\em 
Let $E= E(C(B)_A;V)$. If $V$ is flat as a left $A$-module or $C$ is
coflat as a right $B$-comodule, 
then 
the right $B$-comodules ${\rm Hom}^{B-}(B,E)$ and   
${\rm Hom}^\psi (C,V)$ are isomorphic to each other. }\vspace{.5\baselineskip}  
  
\note{\proof The isomorphism and its inverse are: 
$  
\theta: {\rm Hom}^\psi(C,V)\to {\rm Hom}^{B-}(B,E)$,  
$\theta(\phi)(b) = c\sw 1\otimes_A\phi(c\sw 2)$, where
$c\in\pi_C^{-1}(b)$, $\theta^{-1}(s) =
{}_V\mu\circ(\check{\tau}\otimes_A  V)\circ(C\tens  
s\circ\pi_C)\circ\Delta$. These maps are well-defined by
(\ref{con.phi}*) and the faithful flatness of $V$. 
\endproof \vspace{.5\baselineskip}}  
\noindent{\normalsize\sc Proposition~\ref{crit.ff}*} {\em   
If an  $A$-Galois coextension $C(B)_A$ admits a counital $B$-bicomodule map  
$s :B\to C$ then $C$ is faithfully coflat as a left $B$-comodule.}  
\vspace{.5\baselineskip} 
 
\noindent{\normalsize\sc Proposition~\ref{pro.sec.cleft}*} {\em    
An $A$-Galois coextension $C(B)_A$ is cocleft if and only if  
there exists    
a cross-section $s\in{\rm Hom}^{B-}(B,E)$ such that $\hat{s} 
:= \mu_A\circ(s\tens A):   
B\tens A\to C$ is a bijection.}\vspace{.5\baselineskip}  
   
\noindent\proof Given $s$ with bijective  $\hat s$, the  
cocleaving map and its convolution inverse are   
$ 
\Phi = (\eps\otimes A)\circ\hat{s}^{-1}$, and $\Phi^{-1} =  
\check{\tau}\circ (C\otimes s\circ\pi_C)\circ\Delta.  
$ \endproof 
  
An  
$A$-Galois coextension $C(B)_A$ is cocleft if and only if  
$C\cong    
B\tens A$ in ${}^B\CM_A$.
\vspace{.3\baselineskip} 
 
\noindent{\normalsize\sc Definition~\ref{vr}*} {\em    
Let $C(B)_A$ be an algebra Galois coextension. A right $B$-comodule 
$\bar{E}$ is called  
a {\rm right comodule associated to $C(B)_A$} iff there exists a   
right $A$-module $V$ such that $\bar{E} = (V\otimes C)^0$, where   
$V\otimes C$ is an $(A,C)_\psi$-module of the canonical
entwining structure as in  
Corollary~\ref{MC1}(1). In this case $\bar{E}$ is denoted by
$\bar{E}(V;C(B)_A)$.}\vspace{.5\baselineskip} 
 
\noindent{\normalsize\sc Proposition~\ref{triv.r}*} {\em   
If $C(B)_A$ is a cocleft coextension then any $\bar{E}(V;C(B)_A)$ 
is isomorphic to $V\otimes B$ as a right 
$B$-comodule.} \vspace{.5\baselineskip}  
 
\proof The isomorphism and its inverse are:  
$$  
V\otimes B \to \bar{E}, \qquad v\otimes b\mapsto \pi_{V\otimes 
C}(v\cdot\Phi(c\sw 1)\tens c\sw 2), \quad c\in \pi_C^{-1}(b), 
$$  
$$  
\bar{E}\to V\otimes B, \qquad  x\mapsto \sum_i 
v^i\cdot\Phi^{-1}(c^i\sw 1)\otimes \pi_C(c^i\sw 2), \quad \sum_i v^i\otimes  
c^i \in \pi_{V\otimes C}^{-1}(x),
$$     
where $\Phi$ is a cocleaving map. \endproof   
 
\noindent{\normalsize\sc Definition~\ref{sec.r}*  
{\em Let $\bar{E}$ be a right comodule associated to $C(B)_A$. 
Any right $B$-comodule map $s: B\to\bar{E}$ is called a {\rm 
cross-section}   
of $\bar{E}$.}}\vspace{.5\baselineskip}

The space ${\rm Hom}^{-B}(B,\bar{E})$ of cross-sections of $\bar{E}$ has a   
natural left $B$-comodule structure. Let ${\rm Hom}_{-A}(C,V)$ denote  
the space of right $A$-module maps  $\phi: C\to V$. 
Since $\mu_C$ is left-colinear over $B$ the space   
${\rm Hom}_{-A}(C,V)$ is a left $B$-comodule.\vspace{.5\baselineskip} 
  
\noindent{\normalsize\sc Theorem~\ref{cross.r}*} {\em   
Let $\bar{E}= \bar{E}(V;C(B)_A)$.
If $C$ is faithfully coflat as a left $B$-comodule then   
${\rm Hom}_{-A}(C,V)\cong {\rm Hom}^{-B}(B,\bar{E})$ as left   
$B$-comodules. }\vspace{.5\baselineskip}  
 
\noindent{\normalsize\sc Example~\ref{ex.psi-1}*} {\em   
 For a Hopf-Galois coextension $C(B)_H$ the canonical entwining map
$\psi: c\otimes h\mapsto h\sw 1\otimes
 c\cdot h\sw 2$ 
is bijective if and only if the antipode in $H$ is 
bijective. }\vspace{.5\baselineskip}  
 
\noindent{\normalsize\sc Lemma~\ref{lem.psi-1}*} {\em 
Let $C(B)_A$ be an $A$-Galois coextension with the bijective canonical 
entwining map $\psi$. Then: 
 
(1) $C$ is a left $A$-module with the action ${}_C\mu = 
(C\otimes\eps\circ\mu_C)\circ(\Delta\otimes A)\circ \psi^{-1}$. 
 
(2) The canonical right $C$-comodule left $A$-module map 
$cocan_L:A\otimes C\to C\square_B C$, $a\otimes c\mapsto a\cdot c\sw 1\otimes 
c\sw 2$ is bijective. 
 
(3) The coalgebra $B$ is isomorphic to $\bar{B} := C/\bar{I}_C$, where 
$\bar{I}_C := \span\{ a\cdot c - \eps(a\cdot c\sw 1)c\sw 2\; | \; 
\forall a\in A, c\in C\}$.}\vspace{.5\baselineskip} 
 
\proof Assertion (1) can be proven by direct calculations which, in
particular,  use 
the equation $c\cdot a = \eps(c\sw 1\cdot a_\alpha)c\sw 2^\alpha$, 
relating $\mu_C$ with $\psi$. To prove (2) one
directly verifies that $cocan_L  = cocan\circ\psi^{-1}$. To prove (3)
one first defines the map  
$
\bar{\iota}_C : 
A\otimes C \to \bar{I}_C$, 
$ a\otimes c\mapsto a\cdot c - \eps(a\cdot c\sw 1) c\sw 2$.  
An easy calculation reveals that  $\bar{\iota}_C =
-\iota_C\circ\psi^{-1}$, where $\iota_C:C\tens A \to I_C$,
$c\tens a   
\mapsto c\cdot a - c\sw 0 \eps(c\sw 1\cdot a)$,  
and thus $I_C = \bar{I}_C$, i.e. $B=\bar{B}$. \endproof 

If $\psi$ is bijective then $C\in{}_A^C\M(\psi^{-1})$. Therefore  one
can  define a functor $_A^C\M(\psi^{-1})\to  
{}^B\M$, $M\mapsto {}^0M$, where $^0M := M/\bar{I}_M$,  
$$ 
\bar{I}_M := \span\{a\cdot m - \eps(a\cdot m\sw{-1})m\sw 0 \; |\;
\forall m\in M,  
a\in A\}. 
$$ 
If $V$ is a left $A$-module then $C\otimes V$ is an object in
$_A^C\M(\psi^{-1})$  
where the left action is given by $_{C\otimes V}\mu = (C\otimes {}_V\mu) 
\circ(\psi^{-1}\otimes V)$, and the coaction is ${}_{C\otimes V}\Delta(  
c\otimes v) = \Delta\otimes V$. 
\vspace{.5\baselineskip} 

\noindent{\normalsize\sc Proposition~\ref{relation}*} {\em 
Let $C(B)_A$ be an $A$-Galois coextension with the bijective  canonical 
entwining map $\psi$. Then: 
 
(1) For any left $A$-module $V$,  the left $B$-comodules $^0(C\otimes 
V)$ and $C\otimes_AV$ are isomorphic to each other. 
 
(2) For any right $A$-module $V$, the right $B$-comodules $(V\otimes 
C)^0$ and $V\otimes_A C$ are isomorphic to each
other.}\vspace{.5\baselineskip}  
 
\proof (1) Consider left $B$-comodule maps  
$\bar{\iota}_{C\otimes V}: A\otimes 
C\otimes V \to C\otimes V$, $a\otimes c\otimes v \mapsto 
a\cdot(c\otimes v) - \eps(a\cdot c\sw 1)c\sw 2\otimes v$,  and
$\kappa : C\otimes A\otimes V 
\to C\otimes V$, $c\otimes a\otimes v \mapsto
c\cdot a\otimes v -   
c\otimes a\cdot v$. An easy calculation 
shows that $\bar{\iota}_{C\otimes V} = \kappa\circ (\psi^{-1}\otimes V)$.
Therefore we have a
commutative diagram of $B$-comodule maps with exact rows:
$$  
\begin{CD}  
A\otimes C\otimes V  @>{\bar{\iota}_{C\otimes V}}>>   C\otimes V  
  @>>>{}^0(C\otimes V)@>>> 0 \\  
  @VV{\psi^{-1}\otimes V}V   @VV=V @VVV\\  
C\otimes A\otimes V  @>\kappa>>   C\otimes V  
  @>>>C\otimes_A V@>>> 0.
\end{CD}  
$$ 
Thus we conclude that  ${}^0(C\otimes V)\cong C\otimes_A V$ as left
$B$-comodules. 

(2) Follows from (1) by the left-right symmetry. \endproof

\section*{\normalsize\sc\centering  
Acknowledgements}\vspace{-.6\baselineskip}  
 
I would like to thank Piotr Hajac and Shahn Majid for many 
interesting discussions and valuable comments. I am also grateful to
the referee for numerous useful remarks, and for suggesting the present 
proofs of Theorem~\ref{C-aut}
and Example~\ref{ex.psi-1}(2) in particular. 
Research supported by the EPSRC grant GR/K02244 and by the Lloyd's of
London Tercentenary Foundation.


\baselineskip 6pt  
\footnotesize   
  
\begin{thebibliography}{99}     

\note{\vspace*{-1.5mm}\bibitem{AndDev:ext}     
N.~Andruskiewitsch and J.~Devoto.     
\newblock Extensions of Hopf algebras.
\newblock {\em Algebra i Analiz}, 7:22--61, 1995.}
  
\note{\vspace*{-1.5mm}\bibitem{Bou:com}  
N.~Bourbaki.  
\newblock {\em Commutative Algebra}  
\newblock  Addison-Wesley, Reading, 1972.  }
  
\vspace*{-1.5mm}\bibitem{Brz:tra}     
T.~Brzezi\'nski.     
\newblock Translation map in quantum principal bundles.     
\newblock {\em J.~Geom.\ Phys.}, 20:349--370, 1996. 
     
\note{\vspace*{-1.5mm}\bibitem{Brz:ind}     
T.~Brzezi\'nski.     
\newblock Induction functors, Frobenius properties and Maschke-type
theorems for entwined modules.
\newblock {\em In preparation.}, 1998.}

\vspace*{-1.5mm}\bibitem{BrzHaj:coa}     
T.~Brzezi\'nski and P.M.~Hajac.     
\newblock Coalgebra extensions and algebra coextensions of Galois type.     
\newblock {\em Preprint},  q-alg/9708010, 1997.  {\em Commun. Algebra}
to appear.   
     
\vspace*{-1.5mm}\bibitem{BrzMa:gau}     
T.~Brzezi\'nski and S.~Majid.     
\newblock Quantum group gauge theory on quantum spaces.     
\newblock {\em Commun.\ Math.\ Phys.}, 157:591--638, 1993.     
\newblock Erratum 167:235, 1995. 
     
\vspace*{-1.5mm}\bibitem{BrzMa:coa}     
T.~Brzezi\'nski and S.~Majid.     
\newblock Coalgebra bundles.     
\newblock {\it Commun.\ Math.\ Phys.}, 191, 467--492, 1998.     
     
\vspace*{-1.5mm}\bibitem{BrzMa:dif}     
T.~Brzezi\'nski and S.~Majid.     
\newblock Quantum differentials and the q-monopole revisited.     
\newblock {\em Preprint}  q-alg/9706021, 1997. {\em Acta
Appl. Math.} to appear.    
 
\vspace*{-1.5mm}\bibitem{BrzMa:geo}     
T.~Brzezi\'nski and S.~Majid.     
\newblock Quantum geometry of algebra factorisations and coalgebra bundles. 
\newblock {\em Preprint} math.QA/9808067, 1998. 
   
\note{\vspace*{-1.5mm}\bibitem{BudKon:pri}     
R.J. Budzy\'nski and W. Kondracki.  
\newblock Quantum principal fibre bundles: topological aspects  
\newblock {\em Rep. Math. Phys.} 37:365--385, 1996.}  
  
\note{\vspace*{-1.5mm}\bibitem{CaeMil:Doi}
S.~Caeneepel, G.~Militiaru and S.~Zhu.
\newblock Doi-Hopf modules, Yetter-Drinfel'd modules and Frobenius type
properties. 
\newblock {\em Trans. Amer. Math. Soc.}, 349:4311--4342, 1997.}

\note{\vspace*{-1.5mm}\bibitem{CaeMil:Mas}
S.~Caeneepel, G.~Militiaru and S.~Zhu.
\newblock A Maschke-type Theorem for Doi-Hopf modules and Applications.
\newblock {\em J. Algebra}, 187:388-412, 1997.}

\vspace*{-1.5mm}\bibitem{CaeRai:ind}
S.~Caeneepel, S.~Raianu.
\newblock Induction functors for the Doi-Koppinen unified Hopf modules.
\newblock [in:] {\em Abelian Groups and Modules}, A. Facchini and C.
Menini (eds.), Kluwer Academic Publishers, Dordrecht 1995, pp. 73--94.


\vspace*{-1.5mm}\bibitem{ChaSwe:hop}   
S.U.~Chase and M.E.~Sweedler.   
\newblock {\em Hopf Algebras and Galois Theory.}   
\newblock Springer-Verlag, Berlin-Heidelberg-New York, 1969.    

\note{\vspace*{-1.5mm}\bibitem{DabHaj:exp}     
L.~D\c abrowski, P.M.~Hajac and P.~Siniscalco.
\newblock Explicit Hopf-Galois Description of $SL_{e^{2\pi
i/3}}(2)$-induced Frobenius homomorphisms.
\newblock {\em Preprint}, DAMTP/97-96, SISSA 43/97/FM, q-alg/9708031.}
   
\note{\vspace*{-1.5mm}\bibitem{DijKor:hom}     
M.~S.~Dijkhuizen and T.~Koornwinder.     
\newblock Quantum homogeneous spaces, duality and quantum 2-spheres.     
\newblock {\em Geom.\ Dedicata}, 52:291--315, 1994.}     

\vspace*{-1.5mm}\bibitem{Doi:rel}     
Y.~Doi.   
\newblock On the structure of relative Hopf modules.
\newblock {\em Commun. Alg.}, 11:243--255, 1983.   
   
\vspace*{-1.5mm}\bibitem{Doi:uni}
Y.~Doi.   
\newblock Unifying Hopf modules.
\newblock {\em J. Algebra}, 153:373--385, 1992.   

  
\vspace*{-1.5mm}\bibitem{DoiTak:cle}     
Y.~Doi and M.~Takeuchi.   
\newblock Cleft comodule algebras for a bialgebra.   
\newblock {\em Commun. Alg.}, 14:801--817, 1986.  
   
\vspace*{-1.5mm}\bibitem{DoiTak:miy}     
Y.~Doi and M.~Takeuchi.   
\newblock Hopf-Galois extensions of algebras, the Miyashita-Ulbrich  
action, and Azumaya algebras.  
\newblock {\em J. Alg.}, 121:488--516, 1989.   
  
\vspace*{-1.5mm}\bibitem{Dur:geo}   
M.~Durdevi\'c.   
\newblock Geometry of quantum principal bundles I, II,    
\newblock {\em Commun. Math. Phys.}, 175:427--521, 1996;    
\newblock {\em Rev. Math. Phys.}, 9:531--607, 1997.   
   
\vspace*{-1.5mm}\bibitem{Gug:ext}    
V.K.A.M.~Gugenheim.    
\newblock On extensions of algebras, co-algebras and Hopf algebras, I.    
\newblock {\em Amer. J. Math.}, 84:349--382, 1962.    
     
\vspace*{-1.5mm}\bibitem{Haj:str} P.~M.~Hajac.     
Strong connections on quantum principal bundles.     
{\it Commun.\ Math.\ Phys.} 182:579--617, 1996.     
     
\vspace*{-1.5mm}\bibitem{Hus:fib}     
D.~Husemoller.  
\newblock {\it Fibre Bundles.}     
\newblock 3rd edition, Springer-Verlag, New York, 1994. 

\note{\vspace*{-1.5mm}\bibitem{Kar:hom}  
M.~Karoubi.  
\newblock Homologie cyclique et K-th\'eorie.  
\newblock {\em Ast\'erisque}, 149, 1987. } 
  
\vspace*{-1.5mm}\bibitem{KreTak:hop}     
H.F.~Kreimer  and M.~Takeuchi.   
\newblock Hopf algebras and Galois extensions of an algebra.   
\newblock {\em Indiana Univ. Math. J.}, 30:675--691, 1981.   
     
\vspace*{-1.5mm}\bibitem{Kop:var}
M.~Koppinen.
\newblock Variations on the smash product with applications to
group-graded rings. 
\newblock {\em J. Pure Appl. Alg.}, 104:61--80, 1994. 

\vspace*{-1.5mm}\bibitem{Maj:phy}     
S.~Majid.     
Physics for algebraists: Non-commutative and non-cocommutative     
Hopf algebras by a bicrossproduct construction.     
{\em J.~Algebra} 130:17--64, 1990.   
     
\vspace*{-1.5mm}\bibitem{Maj:rie}     
S.~Majid.     
\newblock Quantum and braided group Riemannian geometry.   
\newblock {\em Preprint} DAMTP/97-73, q-alg/9709025.   
   
\note{\vspace*{-1.5mm}\bibitem{Maj:adv}     
S.~Majid.     
\newblock Advances in quantum and braided geometry.     
\newblock To appear in ed. V.~Dobrev, Proc. XXI ICGTMP (Quantum groups     
volume), Goslar 1996, Heron Press, Sofia. (q-alg/9610003) }    
 
\vspace*{-1.5mm}\bibitem{MasDoi:gen}     
A.~Masuoka and Y.~Doi.   
\newblock Generalization of left comodule algebras.   
\newblock {\em Commun. Alg.}, 20:3703--3721, 1992. 
    
\vspace*{-1.5mm}\bibitem{MilMoo:str}    
J.W.~Milnor and J.C.~Moore.    
\newblock On the structure of Hopf algebras.    
\newblock {\em Ann. Math.}, 81:211--264, 1965.    
    
\vspace*{-1.5mm}\bibitem{Mon:hop}     
S.~Montgomery.     
\newblock {\em Hopf Algebras and Their Actions on Rings}.     
\newblock CBMS Lectures vol. 82, AMS, Providence, RI, 1993.    
   
   
\vspace*{-1.5mm}\bibitem{Pfl:fib}   
M.J.~Pflaum.   
\newblock Quantum groups on fibre bundles.   
\newblock {\em Commun. Math. Phys.} 166:279--315, 1994.   
    
\note{\vspace*{-1.5mm}\bibitem{Pod:sph}     
P.~Podle\'s.     
\newblock Quantum spheres.     
\newblock {\em Lett.\ Math.\ Phys.}, 14:193--202, 1987. }    
     
     
\vspace*{-1.5mm}\bibitem{RadTow:yet}    
D.E.~Radford and J.~Towber.    
\newblock Yetter-Drinfeld categories associated to an arbitrary algebra.
\newblock {\em J. Pure Appl. Algebra.}, 87:259--279, 1993.  

\note{\vspace*{-1.5mm}\bibitem{Sch:bia}     
P. Schauenburg.     
\newblock A bialgebra that admits a Hopf-Galois extension is a Hopf
algebra. 
\newblock {\em Proc. Amer. Math. Soc.} 125:83--85, 1997. }    
 

\vspace*{-1.5mm}\bibitem{Sch:pri}     
H.-J. Schneider.     
\newblock Principal homogeneous spaces for arbitrary Hopf algebras.     
\newblock {\em Isr.\ J.\ Math.}, 72:167--195, 1990.     
     
\vspace*{-1.5mm}\bibitem{Sch:rep}     
H.-J.~Schneider.     
\newblock Representation theory of Hopf-Galois extensions.     
\newblock {\em Isr.\ J.\ Math.}, 72:196--231, 1990.     
     
\vspace*{-1.5mm}\bibitem{Sch:nor}     
H.-J. Schneider.     
\newblock Normal basis and transitivity of crossed products for Hopf     
algebras.     
\newblock {\em J.~Algebra}, 152:289--312, 1992.     
     
\note{\vspace*{-1.5mm}\bibitem{Wat:aff}     
W.C.~Waterhouse.  
\newblock {\it Introduction to Affine Group Schemes}     
\newblock Springer-Verlag, New York, Heidelberg, Berlin, 1979. } 
  
\note{\vspace*{-1.5mm}\bibitem{sch}     
H.-J.~Schneider.      
Some remarks on exact sequences of quantum groups.      
{\it Commun.\ Alg.}, 21:3337--3357 (9), 1993}     
     
\vspace*{-1.5mm}\bibitem{swe}     
M.~E.~Sweedler.       
\newblock {\it Hopf Algebras}     
\newblock W.A.Benjamin, Inc., New York, 1969    
     
\note{\vspace*{-1.5mm}\bibitem{wor} S.~L.~Woronowicz.     
Differential calculus on compact matrix pseudogroups (quantum groups).      
{\it Commun.\ Math.\ Phys.} 122:125--170, 1989. }    
     
\vspace*{-1.5mm}\bibitem{Yet:rep}    
D.N.~Yetter.    
\newblock Quantum groups and representations of monoidal categories.
\newblock {\em Math. Proc. Camb. Phil. Soc.}, 108:261--290, 1990. 

\end{thebibliography}
\end{document}